\documentclass[12pt,a4paper, draft]{article}
\usepackage[dvipdfm, margin=2cm]{geometry}
\usepackage[tbtags]{amsmath}%
\usepackage{amsthm}
\usepackage{amssymb}%
\usepackage{amsfonts}
\usepackage{mathrsfs}
\usepackage{color}
\usepackage{mathtools}
\usepackage{verbatim}
\allowdisplaybreaks[4]
\numberwithin{equation}{section}
\newtheorem{theorem}{Theorem}[section]
\newtheorem{proposition}[theorem]{Proposition}
\newtheorem{corollary}[theorem]{Corollary}
\newtheorem{lemma}[theorem]{Lemma}
\theoremstyle{definition}
\newtheorem{remark}[theorem]{Remark}

\newtheorem{definition}[theorem]{Definition}
\newtheorem{notation}[theorem]{Notation}
\DeclareMathOperator{\ad}{ad}

\DeclareMathOperator{\Rad}{Rad}
\DeclareMathOperator{\Res}{Res}

\newcommand{\BF}{\mathbb{F}}
\newcommand{\BN}{\mathbb{N}}
\newcommand{\BZ}{\mathbb{Z}}

\newcommand{\BC}{\mathbb{C}}

\newcommand{\spanf}{\mathrm{span}_\BF}
\newcommand{\Vir}{\mathcal{V}\mathrm{ir}}
\newcommand{\fg}{\mathfrak{g}}
\newcommand{\fh}{\mathfrak{h}}
\newcommand{\D}{\mathcal{D}}

\newcommand{\B}{\mathcal{B}}

\newcommand{\bc}{\mathbf{c}}

\newcommand{\1}{\mathbf{1}}

\newcommand{\?}{\hspace*{1.346em}}
\makeatletter
\newcommand{\leqnomode}{\tagsleft@true}
\newcommand{\reqnomode}{\tagsleft@false}
\makeatother
\makeatletter
\newtoks\@enLab  
\def\@enQmark{?}
\def\@enLabel#1#2{%
	\edef\@enThe{\noexpand#1{\@enumctr}}%
	\@enLab\expandafter{\the\@enLab\csname the\@enumctr\endcsname}%
	\@enloop}
\def\@enSpace{\afterassignment\@enSp@ce\let\@tempa= }
\def\@enSp@ce{\@enLab\expandafter{\the\@enLab\space}\@enloop}
\def\@enGroup#1{\@enLab\expandafter{\the\@enLab{#1}}\@enloop}
\def\@enOther#1{\@enLab\expandafter{\the\@enLab#1}\@enloop}
\def\@enloop{\futurelet\@entemp\@enloop@}
\def\@enloop@{%
	\ifx A\@entemp         \def\@tempa{\@enLabel\Alph  }\else
	\ifx a\@entemp         \def\@tempa{\@enLabel\alph  }\else
	\ifx i\@entemp         \def\@tempa{\@enLabel\roman }\else
	\ifx I\@entemp         \def\@tempa{\@enLabel\Roman }\else
	\ifx 1\@entemp         \def\@tempa{\@enLabel\arabic}\else
	\ifx \@sptoken\@entemp \let\@tempa\@enSpace         \else
	\ifx \bgroup\@entemp   \let\@tempa\@enGroup         \else
	\ifx \@enum@\@entemp   \let\@tempa\@gobble          \else
	\let\@tempa\@enOther
	\fi\fi\fi\fi\fi\fi\fi\fi
	\@tempa}
\newlength{\@sep} \newlength{\@@sep}
\providecommand{\sfbc}{\rmfamily\upshape}
\providecommand{\sfn}{\rmfamily\upshape}
\def\@enfont{\ifnum \@enumdepth >1\let\@nxt\sfn \else\let\@nxt\sfbc \fi\@nxt}
\def\enumerate{%
	\ifnum \@enumdepth >3 \@toodeep\else
	\advance\@enumdepth \@ne
	\edef\@enumctr{enum\romannumeral\the\@enumdepth}\fi
	\@ifnextchar[{\@@enum@}{\@enum@}}
\def\@@enum@[#1]{%
	\@enLab{}\let\@enThe\@enQmark
	\@enloop#1\@enum@
	\ifx\@enThe\@enQmark\@warning{The counter will not be printed.%
		^^J\space\@spaces\@spaces\@spaces The label is: \the\@enLab}\fi
	\expandafter\edef\csname label\@enumctr\endcsname{\the\@enLab}%
	\expandafter\let\csname the\@enumctr\endcsname\@enThe
	\csname c@\@enumctr\endcsname7
	\expandafter\settowidth
	\csname leftmargin\romannumeral\@enumdepth\endcsname
	{\the\@enLab\hskip\labelsep}%
	\@enum@}
\def\@enum@{\list{{\@enfont\csname label\@enumctr\endcsname}}%
	{\usecounter{\@enumctr}\def\makelabel##1{\hss\llap{##1}}%
		\ifnum \@enumdepth>1\setlength{\topsep}{\@@sep}\else
		\setlength{\topsep}{\@sep}\fi
		\ifnum \@enumdepth>1\setlength{\itemsep}{0pt plus1pt minus1pt}%
		\else \setlength{\itemsep}{\@@sep}\fi
		\setlength{\parsep}{0pt plus1pt minus1pt}%
		\setlength{\parskip}{0pt plus1pt minus1pt}
}}

\def\endenumerate{\par\ifnum \@enumdepth >1\addvspace{\@@sep}\else
	\addvspace{\@sep}\fi \endlist}
\makeatother                                                            %

\begin{document}
	\title{Modular vertex algebras associated with the centrally extended planar Galilean conformal algebra and their representations}
	\author{
		Jinxin Hu$^{a}$ \quad and \quad Hongju Zhao$^{b}$\\
		\date{ {$^{a}$Department of Mathematics, Soochow University, Suzhou 215006, China}\\
			{$^{b}$School of Mathematical Sciences, Harbin Normal University, Harbin 150025, China}
	}}
	\maketitle

\begin{abstract}
We study the vertex algebra $V_{\mathcal{G}}(c,0)$ associated with $\mathcal{G}$ over an algebraically closed field $\BF$ of characteristic $p>2$, 
where $\mathcal{G}$ is a modular analogue  of the universal central extension of the planar Galilean conformal algebra
over $\BC$ and $c=(c_1,c_2,c_3)\in\BF^3$.
We first classify the irreducible $\BN$-graded modules of $V_{\mathcal{G}}(c,0)$ and of its graded quotient $V_{\mathcal{G}}^0(c,0)$;
in particular, $V_{\mathcal{G}}^0(c,0)$ has exactly $p^2$ irreducible $\BN$-graded modules up to equivalence.
We then show that the space of invariant bilinear forms on $V_{\mathcal{G}}(c,0)$ is one-dimensional when $c_2=0$ and zero otherwise.
For $c_2=0$, 
the quotient by the radical of a nonzero invariant bilinear form is the simple graded quotient,
which is isomorphic to the simple Virasoro vertex algebra of central charge $c_1$ when $c_3=0$, and to the tensor product of the simple graded quotient of the rank-one Heisenberg vertex algebra at level $c_3$ and the simple Virasoro vertex
algebra of central charge $c_1-1$ when $c_3\ne0$.
Finally, using these identifications, 
we classify the irreducible $\BN$-graded modules of the simple graded quotient when the corresponding Virasoro vertex algebra has central charge $\frac{1}{2}$.	
\end{abstract}

\section{Introduction}
Vertex algebras were introduced by Borcherds in 1986 (see \cite{B1}), with axioms formulated over an arbitrary field.
Among the most important examples are those constructed from infinite-dimensional Lie algebras, such as affine Lie algebras, Heisenberg algebras, and the Virasoro algebra (see \cite{FZ}).
The representation theory of these vertex algebras in prime characteristic has received considerable attention (see, e.g., \cite{ATV,AW,DR,DR2015,DLR1,JLM,LM}).

Several general theories familiar from characteristic zero have also
been developed for vertex algebras over fields of prime characteristic, 
including those concerning Zhu algebras (see \cite{DR,JLM}), local systems of vertex operators (see \cite{L1996,LM}), and invariant bilinear forms (see \cite{LM2}).
On the other hand, concrete examples exhibit interesting new features in prime characteristic.
For instance, unlike in characteristic zero, universal Heisenberg vertex algebras at nonzero levels are not simple.
Li and Mu constructed a family of simple quotients
and showed that each quotient has a unique irreducible module up
to isomorphism and that every module for the quotient is completely
reducible (see \cite{LM}).
Moreover, quotients of the universal Virasoro and affine vertex
algebras by ideals associated with the $p$-centers have been
studied. 
The Zhu algebras of certain graded quotients were determined and used to classify the irreducible $\BN$-graded modules for those quotients (see \cite{JLM}).

The planar Galilean conformal algebra is the
infinite-dimensional Galilean conformal algebra in
$(2+1)$-dimensional spacetime arising from Bagchi and
Gopakumar's study of a nonrelativistic limit of the AdS/CFT
correspondence (see \cite{BG}).
Its universal central extension $\mathcal{G}_{\BC}$ over $\BC$
was determined in \cite{GLP}.
The Lie algebra $\mathcal{G}_{\BC}$ contains a twisted Heisenberg--Virasoro subalgebra and two subalgebras of $\mathcal{W}(2,2)$-type (see \cite{G2026}).
Motivated by these connections, we consider the
modular analogue $\mathcal{G}$ of $\mathcal{G}_{\BC}$ over an
algebraically closed field $\BF$ of characteristic $p>2$ and
study the associated vertex algebras.

The Lie algebra $\mathcal{G}$, which we call the
	modular centrally extended planar Galilean conformal algebra,
	has basis
\begin{equation*}
	\{L_n,H_n,I_n,J_n\mid n\in\BZ\}
	\cup\{\bc_1,\bc_2,\bc_3\},
\end{equation*}
where $\bc_1,\bc_2,\bc_3$ are central,
and the nonzero Lie brackets, up to skew-symmetry, are given by
\eqref{eq:Lm-Ln-brac}--\eqref{eq:H-I;H-J}.

Let $\mathfrak{p}$ be the subalgebra
of $\mathcal{G}$ spanned by $L_n,I_n,J_n$ for $n\ge -1$ and $H_m$ for $m\ge 0$.
For $c=(c_1,c_2,c_3)\in\BF^3$, let $\BF_c$ be the
one-dimensional
$\bigl(\mathfrak p\oplus\bigoplus_{i=1}^3\BF\bc_i\bigr)$-module
on which $\mathfrak p$ acts trivially and $\bc_i$
acts as $c_i$ for $i=1,2,3$.
Set
\begin{equation*}
	V_{\mathcal{G}}(c,0)
	:=U(\mathcal{G})
	\otimes_{U(\mathfrak{p}\oplus\bigoplus_{i=1}^3\BF\bc_i)}
	\BF_c.
\end{equation*}
Theorem~\ref{prop:Galilean-VA} gives $V_{\mathcal{G}}(c,0)$
a vertex algebra structure uniquely determined by the vacuum vector
$\1=1\otimes 1$ and the vertex operators in \eqref{eq:ver-ope-2}.
Moreover, its modules are in one-to-one correspondence with
locally truncated $\mathcal{G}$-modules of central charge $c$ (see Proposition~\ref{prop:mod-equiv}).

We then show that $\mathcal{G}$ carries a restricted Lie algebra
structure (see Lemma~\ref{lem:rest-stru-on-G}) and use the
$p$-center of $U(\mathcal{G})$ to construct a family of quotient
vertex algebras of $V_{\mathcal{G}}(c,0)$.
More precisely, let $\fg_{\mathcal{G}}$ be the Lie subalgebra spanned by $L_0,H_0,I_0,J_0$.
For a linear functional
$\chi:\fg_{\mathcal{G}}\to\BF$,
let $I^{\chi}$ be the $\mathcal{G}$-submodule of $V_{\mathcal{G}}(c,0)$ generated by the elements in \eqref{eq:gen-I-chi}.
It is proved that $I^{\chi}$ is an ideal of $V_{\mathcal{G}}(c,0)$ (see Theorem~\ref{th:I-mu-iedal}).
Thus, we obtain the quotient vertex algebra
\begin{equation*}
V^{\chi}_{\mathcal{G}}(c,0):=V_{\mathcal{G}}(c,0)/I^{\chi}.
\end{equation*}
Furthermore,
the $V^\chi_{\mathcal{G}}(c,0)$-modules
are in one-to-one correspondence with
locally truncated $\mathcal{G}$-modules of central charge $c$
satisfying the relations in \eqref{eq:L-2}--\eqref{eq:H-1}
(see Proposition~\ref{prop:quot-mod-equiv}).

For $\chi=0$, the quotient $V^0_{\mathcal{G}}(c,0)$ is
$\BN$-graded, as is $V_{\mathcal{G}}(c,0)$.
We prove that the Zhu algebras of $V_{\mathcal{G}}(c,0)$ and
$V^0_{\mathcal{G}}(c,0)$ are isomorphic to $U(\fg_{\mathcal{G}})$
and $\mathfrak u(\fg_{\mathcal{G}})$, respectively
(see Theorem~\ref{th:A(V)} and
Theorem~\ref{th:iso-u(g)-V0}),
where $\mathfrak u(\fg_{\mathcal{G}})$ is the restricted
enveloping algebra of $\fg_{\mathcal{G}}$.
We then classify the irreducible $\fg_{\mathcal{G}}$-modules and determine which of them factor through
$\mathfrak{u}(\fg_{\mathcal{G}})$
(see Theorem~\ref{th:classification-g-G} and
Corollary~\ref{coro:zero-mode-restricted}).
By Zhu algebra theory (see \cite{Zhu1996}, cf. \cite[Propositions 2.10, 2.11]{JLM}), 
these results give a classification of the irreducible $\BN$-graded modules for	$V_{\mathcal{G}}(c,0)$ and $V^0_{\mathcal{G}}(c,0)$.
In particular, $V^0_{\mathcal{G}}(c,0)$ has exactly $p^2$ irreducible $\BN$-graded modules up to equivalence (see Theorem \ref{th:irr-V0-mod}).

We further show that $V_{\mathcal{G}}(c,0)$ is an $\mathcal{H}$-module
vertex algebra,	
where $\mathcal{H}$ is the Hopf algebra introduced by Li and Mu in \cite{LM2}.
Using this $\mathcal{H}$-module vertex algebra structure, 
we show that the space of invariant bilinear forms on $V_{\mathcal{G}}(c,0)$ is one-dimensional when $c_2=0$ and zero otherwise (see Theorem \ref{th:Galilean-invariant-forms} and Remark \ref{rmk:c2 ne 0}). 
When $c_2=0$, the radical of the unique symmetric invariant bilinear form $(\cdot,\cdot)$ such that $(\1,\1)=1$ is the unique maximal proper graded ideal (see Lemma \ref{lem:rad-max-ideal}). 
We denote the corresponding simple graded quotient by $L_{\mathcal{G}}(c,0)$.

For $c_2=0$, we obtain
\begin{equation*}
	L_{\mathcal{G}}(c,0)\cong
	\begin{cases}
		L_{\Vir}(c_1,0), & \text{if }c_3=0,\\
		L_{\widehat{\fh}}(c_3,0)\otimes L_{\Vir}(c_1-1,0),
		& \text{if }c_3\ne0,
	\end{cases}
\end{equation*}
as $\BN$-graded vertex algebras (see Theorem~\ref{th:c2=0}), where
$L_{\Vir}(c_1,0)$ denotes the simple Virasoro vertex algebra
of central charge $c_1$, and
$L_{\widehat{\fh}}(c_3,0)$ denotes the simple graded quotient of
the rank-one Heisenberg vertex algebra at level $c_3$.
The irreducible $\BN$-graded modules for $L_{\widehat{\fh}}(c_3,0)$ with $c_3\ne0$ and for
$L_{\Vir}(\frac12,0)$ have been classified in \cite{LM} and \cite{DR2015},
respectively.
These works also establish complete reducibility of
$\BN$-graded modules when $p\ne7$.
Consequently, for $c=(\frac12,0,0)$ or
$c=(\frac32,0,c_3)$ with $c_3\ne0$,
$L_{\mathcal{G}}(c,0)$ has exactly three irreducible
$\BN$-graded modules up to equivalence if $p\ne7$
and exactly two if $p=7$.
For these choices of $c$, every $\BN$-graded
$L_{\mathcal{G}}(c,0)$-module is completely reducible
when $p\ne7$ (see Corollary~\ref{coro:c2=0}).

The paper is organized as follows.
In Section~\ref{sec:2}, we recall some basic notions and review the general construction of vertex algebras and their modules.
In Section~\ref{sec:3}, we construct the vertex algebra $V_{\mathcal{G}}(c,0)$ and a family of quotient vertex algebras using the $p$-center of $U(\mathcal{G})$.
In Section~\ref{sec:4}, we determine the Zhu algebras of $V_{\mathcal{G}}(c,0)$ and $V^0_{\mathcal{G}}(c,0)$, and classify their irreducible $\BN$-graded modules.
In Section~\ref{sec:5}, we determine the invariant bilinear forms on $V_{\mathcal{G}}(c,0)$, identify its simple graded quotient when $c_2=0$, and classify the irreducible $\BN$-graded modules of this quotient for certain central charges.

\section{Preliminaries}\label{sec:2}
Throughout this paper, $\BF$ denotes an algebraically closed field of characteristic $p>2$. 
Unless otherwise stated, all vector spaces and algebras are over $\BF$. 
We write $\BC$, $\BZ$, $\BZ_+$, 
and $\BN$ for the sets of complex numbers, 
integers, positive integers, and nonnegative integers, 
respectively.

For any $m\in \BZ$ and $k\in \BN$,
\begin{equation*}
	\binom{m}{k}=\frac{m(m-1)\cdots (m+1-k)}{k!}\in \BZ.
\end{equation*}
Then we shall also view $\binom{m}{k}$ as an element of $\BF$.
We will use {\em Vandermonde's identity}:
\begin{equation}\label{eq:Vand-iden}
	\sum_{i=0}^{k}\binom{m}{i}\binom{n}{k-i}
	=
	\binom{m+n}{k}\quad \text{for }m,n\in\BZ \text{ and }k\in \BN.
\end{equation}

Let $\B$ be the abstract (cocommutative) bialgebra with basis
$\{\D^{(n)}\mid n\in\BN\}$ (see \cite{LM}; cf. \cite{B1}),
where
\begin{align}
	\D^{(m)}\cdot\D^{(n)}&=\binom{m+n}{m}\D^{(m+n)}, \quad \D^{(0)}=1,\label{eq:bi-alg-prod}\\
	\Delta(\D^{(n)})&=\sum_{i=0}^{n}\D^{(n-i)}\otimes \D^{(i)},\quad \varepsilon(\D^{(n)})=\delta_{n,0}\label{eq:bi-alg-coprod}
\end{align}
for $m,n\in \BN$.
Following \cite{DR}, 
set a generating function
\begin{equation*}
	e^{x\D}=\sum_{n\in \BN} x^{n}\D^{(n)}\in \B [[x]].
\end{equation*}
Then \eqref{eq:bi-alg-prod}--\eqref{eq:bi-alg-coprod} can be written as follows:
\begin{equation}
	e^{x\D}e^{z\D}=e^{(x+z)\D},\quad \Delta(e^{x\D})=e^{x\D}\otimes e^{x\D}, \quad \varepsilon(e^{x\D})=1.
\end{equation}

For any vector space $U$, the space $U[[x,x^{-1}]]$ carries a $\B$-module structure with $\D^{(k)}$ for $k\in \BN$ acting as the $k$-th Hasse derivative $\partial_{x}^{(k)}$ (with respect to $x$), which is defined by 
\begin{equation}
	\partial_{x}^{(k)}x^{m}=\binom{m}{k}x^{m-k}\quad \text{for }m\in \BZ.
\end{equation}
Write $e^{z\partial_{x}}:=\sum_{k\in \BN}z^{k}\partial_{x}^{(k)}$.
Then
\begin{equation}\label{eq:taylor}
	e^{z\partial_{x}}f(x)=f(x+z)\quad\text{for }f(x)\in U[[x,x^{-1}]]
\end{equation}
(see, for example, \cite[Remark 2.1]{LM2}).

Let $H$ be a cocommutative bialgebra.
Recall that an {\em $H$-module algebra} is a unital algebra $A$ equipped with an $H$-module structure such that
\begin{equation*}
	h\cdot (ab)=\sum (h_{(1)}\cdot a)(h_{(2)} \cdot b),\quad h\cdot 1=\varepsilon(h)1
\end{equation*}
for $h\in H$ and $a,b\in A$, where $\Delta(h)=\sum h_{(1)}\otimes h_{(2)}$.
An {\em $H$-module Lie algebra} is a Lie algebra $\fg$ equipped with a $H$-module structure such that
\begin{equation}\label{eq:Bialg-mod-com}
	h\cdot [a,b]=\sum [h_{(1)}a,h_{(2)}b]
\end{equation}	
for $h\in H$, $a,b\in \fg$.
Let $A$ be an $H$-module algebra. 
An {\em $(A,H)$-module} (see \cite{LM}) is an $A$-module $W$ equipped with an $H$-module structure such that 
\begin{equation}\label{eq:A-H}
	h(aw)=\sum (h_{(1)}a)(h_{(2)}w)\quad\text{for }h\in H, a\in A, w\in W.
\end{equation}

\begin{lemma}[{\cite[Lemma 3.2]{LM}}]\label{lem:H-mod-Lie-H-mod-alg}
	Let $\fg$ be an $H$-module Lie algebra. 
	Then the universal enveloping algebra $U(\fg)$ is naturally an $H$-module algebra.
	Moreover, in case $H=\B$, for $u\in U(\fg)$, we have
	\begin{equation}
		e^{x\D} u e^{-x\D}=(e^{x\D} u)
	\end{equation}
	on any $(U(\fg), \B)$-module, in particular on $U(\fg)$.
\end{lemma}	

\begin{definition}[cf. \cite{B1,LL}]\label{def:VA}
A {\em vertex algebra} over $\BF$ is a vector space $V$ equipped with a distinguished vector $\1\in V$ (called the {\em vacuum vector}), and a linear map
\begin{align*}
		Y(\cdot,x)\colon V&\to (\operatorname{End} V)[[x,x^{-1}]],\\
		    &v\mapsto Y(v,x)=\sum_{n\in\BZ}v_nx^{-n-1},
\end{align*}
satisfying the following conditions for all $u,v\in V$:
\begin{enumerate}[(1)]
\item (the {\em truncation condition}) $u_nv=0$ for all sufficiently large $n$.
\item (the {\em vacuum property})  $Y(\1,x)=\operatorname{id}_V$.
\item (the {\em creation property}) $Y(v,x)\1\in V[[x]]$ and $\lim\limits_{x\to 0}Y(v,x)\1=v$.
\item (the {\em Jacobi identity})
\begin{align}\label{eq:Jaco-ide}
&\?x_0^{-1}\delta\left(\frac{x_1-x_2}{x_0}\right)
			Y(u,x_1)Y(v,x_2) -x_0^{-1}\delta\left(\frac{x_2-x_1}{-x_0}\right)
			Y(v,x_2)Y(u,x_1)\nonumber\\
			& =x_2^{-1}\delta\left(\frac{x_1-x_0}{x_2}\right)
			Y(Y(u,x_0)v,x_2).
		\end{align}
	\end{enumerate}
\end{definition}
It is known (see \cite[Lemma 2.2]{LM2}; cf. \cite{B1}) that a vertex algebra $V$ carries a $\B$-module structure with 
\begin{equation}\label{eq:B-mod-str}
	\D^{(n)}v=v_{-n-1}\1 \quad\text{ for } n\in\BN.
\end{equation}	

For a vertex algebra $V$, a {\em $V$-module} is a vector space $W$ equipped with a linear map
\begin{align*}
	Y_W(\cdot,x):V&\to (\operatorname{End} W)[[x,x^{-1}]],\\
	&v\mapsto Y_{W}(v,x)=\sum_{n\in\BZ}v_nx^{-n-1},
\end{align*}
satisfying the conditions that $Y_{W}(\1,x)=\operatorname{id}_{W}$, $Y_{W}(v,x)w\in W((x))$ for $v\in V$, $w\in W$, and that the corresponding Jacobi identity \eqref{eq:Jaco-ide} holds.


We end this section with a brief review of the theory of local systems of vertex operators (cf. \cite{L1996}, \cite{LL}).
For any vector space $W$, define
\begin{equation*}
	\mathcal{E}(W):=\operatorname{Hom}(W,W((x))).
\end{equation*}	
Following \cite[Section 3]{L1996}, for $a(x),b(x)\in\mathcal E(W)$ and
$n\in\BZ$, define
\begin{equation*}
	a(x)_n b(x)
	:=
	\operatorname{Res}_{x_1}
	\big(
	(x_1-x)^n a(x_1)b(x)
	-
	(-x+x_1)^n b(x)a(x_1)
	\big).
\end{equation*}
The elements $a(x),b(x)\in\mathcal E(W)$ are said to be {\em mutually
 local} if there exists $k\in\BN$ such that
\begin{equation*}
	(x_1-x_2)^k[a(x_1),b(x_2)]=0.
\end{equation*}
A subset $S\subseteq\mathcal E(W)$ is {\em local} if every pair $a(x),b(x)\in S$ is mutually local.

\begin{theorem}[{\cite[Remark 2.15]{LM}; cf. \cite[Theorem 5.5.18]{LL}}]\label{th:local-set-VA}
	Let $S$ be a local subset of $\mathcal{E}(W)$.
	Then $S$ generates a vertex algebra $\langle S\rangle\subseteq\mathcal E(W)$, and $W$ is a natural faithful $\langle S\rangle$-module.
	The module structure is given by
	$Y_W(a(x),z)=a(z)$ for $a(x)\in\langle S\rangle$.
\end{theorem}

The following construction theorem is the prime-characteristic analogue of \cite[Theorem 5.7.1]{LL}.

\begin{theorem}[{\cite[Theorem 2.14]{LM}}]\label{th:LM-general-VA}
Let $V$ be a $\B$-module with a distinguished vector $\1$ such that $\D^{(n)}\mathbf{1}=0$ for $n\in \BZ_+$ and let $T$ be a subset of $V$ equipped with a map
\begin{equation}\label{eq:def-Y-0}
	Y_0:T\to\operatorname{Hom}(V,V((x))),
	\quad a\mapsto Y_0(a,x)=\sum_{n\in\BZ}a_nx^{-n-1}.
\end{equation}
Assume that the following conditions hold:
\begin{enumerate}[(i)]
	\item For $a\in T$, $Y_{0}(a,x)\mathbf{1}\in V[[x]]$ and $\lim\limits_{x\to 0}Y_{0}(a,x)\mathbf{1}=a$.
	\item For $a,b\in T$, there exists $k\in \BN$ such that
	\begin{equation*}
		(x_1-x_2)^{k}[Y_{0}(a,x_1),Y_0(b,x_2)]=0.
	\end{equation*}
	\item $V$ is linearly spanned by $a_{n_1}^{(1)}\cdots a_{n_r}^{(r)}\mathbf{1}$ for $r\in \BN$, $a^{(1)},\ldots, a^{(r)}\in T$, and $n_1,\ldots,n_r\in \BZ$.
	\item For $a\in T$,
	\begin{equation*}
		e^{z\D}Y_{0}(a,x)e^{-z\D}=e^{z\partial_{x}}Y_{0}(a,x).
	\end{equation*}
\end{enumerate}
Then $Y_{0}$ can be extended uniquely to a linear map $Y:V\to \operatorname{Hom}(V,V((x)))$ such that $(V,Y,\mathbf{1})$ is a vertex algebra.
\end{theorem}

We now give a construction theorem for modules over a vertex algebra
(cf. \cite[Theorem 5.7.13]{LL}).

\begin{theorem}\label{th:gen-con-V-mod}
Let $V$ be a vertex algebra and let $T$ be a generating subset
of $V$. Let $W$ be a vector space equipped with a linear map
\begin{equation}
	Y_W\colon V\to \operatorname{Hom}(W,W((x))),\quad v\mapsto v_W(x).
\end{equation}
Assume that the following conditions hold:
\begin{enumerate}[(i)]
	\item $\1_W(x)=\operatorname{id}_W$.
	\item $(a_nv)_W(x)=a_W(x)_nv_W(x)$ for
	$a\in T$, $v\in V$, and $n\in\BZ$.
	\item For $a,b\in T$, there exists $k\in\BN$ such that
	\begin{equation}\label{eq:2.4.8}
		(x_1-x_2)^k[a_W(x_1),b_W(x_2)]=0.
	\end{equation}
\end{enumerate}
Then $(W,Y_W)$ is a $V$-module.
\end{theorem}

\begin{proof}
Set $V(x)=\{v_W(x)\mid v\in V\}$ and 	$T(x)=\{a_W(x)\mid a\in T\}$.
By (iii), $T(x)$ is a local subset of $\mathcal E(W)$.
Hence Theorem \ref{th:local-set-VA} gives a vertex algebra $\langle T(x)\rangle\subseteq\mathcal E(W)$ generated by $T(x)$,
with $W$ as a natural faithful module.

Since $T$ generates $V$, 
it follows from \cite[Proposition 3.9.3]{LL} that $V$ is spanned by elements of the form  $a^{(1)}_{n_1}\cdots a^{(k)}_{n_k}\1$ with $a^{(i)}\in T$ and $n_i\in \BZ$.
Consequently, (i) and (ii) give $V(x)\subseteq\langle T(x)\rangle$.

We next show that $V(x)$ is a vertex subalgebra containing
$T(x)$, and hence $V(x)=\langle T(x)\rangle$.
Set
\begin{equation*}
	U=\{v\in V\mid
	(v_nw)_W(x)=v_W(x)_nw_W(x)
	\text{ for all }w\in V,\ n\in\BZ\}.
\end{equation*}
Then $U$ is a subspace of $V$, and conditions (i) and (ii)
give $\1\in U$ and $T\subseteq U$.
Since $T$ generates $V$, it suffices to show that
$u_mv\in U$ for $u,v\in U$ and $m\in\BZ$.
For $w\in V$ and $n\in\BZ$, the iterate formula
(see, for example, \cite[(3.1.12)]{LL}) in
$V$ and $\langle T(x)\rangle$, together with $u,v\in U$, gives
\begin{align*}
	&\?((u_{m}v)_{n}w)_{W}(x)\nonumber\\
	&=\sum_{i\in \BN}(-1)^{i}\binom{m}{i}(u_{m-i}v_{n+i}w-(-1)^{m}v_{m+n-i}u_iw)_{W}(x) \nonumber\\
	&=\sum_{i\in \BN}(-1)^{i}\binom{m}{i}(u_{W}(x)_{m-i}v_{W}(x)_{n+i}w_{W}(x)-(-1)^{m}v_W(x)_{m+n-i}u_{W}(x)_iw_{W}(x))\nonumber\\
	&=(u_{W}(x)_{m}v_{W}(x))_{n}w_{W}(x)\nonumber\\
	&=(u_{m}v)_{W}(x)_{n}w_{W}(x).
\end{align*}
Thus $U=V$, so $V(x)$ is a vertex subalgebra containing
$T(x)$ and hence equals $\langle T(x)\rangle$.
Together with (i), this shows that
$Y_W\colon V\to\langle T(x)\rangle$
is a vertex algebra homomorphism.
Hence $(W,Y_W)$ is a $V$-module.
\end{proof}

\section{The vertex algebra $V_{\mathcal{G}}(c,0)$ and its quotients $V_{\mathcal{G}}^{\chi}(c,0)$}\label{sec:3}
In this section, 
we construct and study the vertex algebra $V_{\mathcal{G}}(c,0)$ associated with $\mathcal{G}$,
where $\mathcal{G}$ is the modular analogue over $\BF$ of the universal central extension over $\BC$ of the
planar Galilean conformal algebra
(see \cite{GLP,G2026}).
We also study the quotients $V_{\mathcal{G}}^{\chi}(c,0)$ by ideals associated with the $p$-center of $U(\mathcal{G})$.

\begin{definition}[{cf. \cite[Definition 2.1]{G2026}}]\label{def:UPGCA}
Let $\mathcal{G}$ be the infinite-dimensional Lie algebra over $\BF$ with a basis
\begin{equation*}
\{L_n,H_n,I_n,J_n\mid n\in\BZ\}\cup \{\bc_1,\bc_2,\bc_3\}
\end{equation*}
subject to the nontrivial Lie brackets
\begin{align}
	[L_m,L_n]&=(n-m)L_{m+n}
	+\frac{1}{2}\binom{m+1}{3}\delta_{m+n,0}\bc_1,
	\label{eq:Lm-Ln-brac}\\
	[L_m,H_n]&=nH_{m+n}
	+m(m+1)\delta_{m+n,0}\bc_2,
	\quad
	[H_m,H_n]=m\delta_{m+n,0}\bc_3,
	\label{eq:Hm-Hn-brac}\\
	[L_m,I_n]&=(n-m)I_{m+n},
	\quad
	[L_m,J_n]=(n-m)J_{m+n},
	\label{eq:L-I;L-J}\\
	[H_m,I_n]&=I_{m+n},
	\quad
	[H_m,J_n]=-J_{m+n},
	\label{eq:H-I;H-J}
\end{align}
for $m,n\in\BZ$, 
where $\bc_1,\bc_2,\bc_3$ are central.
\end{definition}

\begin{remark}\label{rem:basis-change}
	The Lie algebra $\mathcal{G}$ is the modular analogue of
	the universal central extension of the planar Galilean
	conformal algebra over $\BC$
	(see \cite{GLP,G2026}), after the change of basis
	\begin{equation*}
		H_n=H_n'+\delta_{n,0}\bc_2
		\quad\text{for }n\in\BZ,
	\end{equation*}
	with all other basis elements unchanged.
	Here $H_n'$ denotes the corresponding generator in
	\cite[Definition~2.1]{G2026}.
	We call $\mathcal{G}$ the
	{\em modular centrally extended planar Galilean conformal algebra}.
\end{remark}

A $\mathcal{G}$-module $W$ is called a
{\em locally truncated module} if for every $w\in W$ and $X=L,I,J,H$,
$X(n)w=0$ for sufficiently large $n$,
where $X(n)$ denotes the operator on $W$
corresponding to $X_n$.
On the other hand,
a $\mathcal{G}$-module $W$ is said to be of
{\em central charge} $c=(c_1,c_2,c_3)\in\BF^3$
if $\bc_i$ acts on $W$ as the scalar $c_i$
for $1\leq i\leq3$.

\begin{lemma}\label{lem:U(G)-B-mod-alg}
	The Lie algebra $\mathcal{G}$ is a $\B$-module Lie algebra
	with the action given by
	\begin{equation}\label{eq:B-mod-acts-G}
		\D^{(k)}(X_m)=(-1)^k\binom{m+1}{k}X_{m-k},\quad
		\D^{(k)}(H_m)=(-1)^k\binom{m}{k}H_{m-k},\quad
		\D^{(k)}(\bc_i)=\delta_{k,0}\bc_i
	\end{equation}
	for $k\in\BN$, $m\in\BZ$, $X=L,I,J$, and $i=1,2,3$.
	Furthermore, $U(\mathcal{G})$ is a $\B$-module algebra.
\end{lemma}

\begin{proof}
A direct calculation shows that $\mathcal{G}$ with the
action \eqref{eq:B-mod-acts-G} is a $\B$-module.
It remains to check that
\begin{equation}\label{eq:Lie-B-com}
	\D^{(k)}[a,b]
	=
	\sum_{i=0}^{k}
	[\D^{(k-i)}a,\D^{(i)}b]
	\quad\text{for }k\in\BN,\ a,b\in\mathcal{G}.
\end{equation}
By bilinearity and skew-symmetry, it suffices to check
\eqref{eq:Lie-B-com} on pairs of basis elements.
The case $k=0$ is clear.
	
Let $k\in\BZ_+$ and $m,n\in\BZ$.
For $a=L_m$ and $b=L_n$, \eqref{eq:Lie-B-com} follows from
\cite[Lemma 3.2]{JLM}, using $-L_m$ as the Virasoro generators.
Recall from \cite[(A.1)]{JLM} that
	\begin{equation}\label{eq:JLM-A1}
		\sum_{i=0}^{k}
		(m-n-k+2i)\binom{m+1}{k-i}\binom{n+1}{i}
		=(m-n)\binom{m+n+1}{k}.
	\end{equation}
For $a=L_m$ and $b=X_n$ with $X=I,J$, 
by \eqref{eq:B-mod-acts-G},
\eqref{eq:L-I;L-J},
and \eqref{eq:JLM-A1}, 
we obtain
\begin{align*}
	\sum_{i=0}^{k}
	[\D^{(k-i)}L_m,\D^{(i)}X_n]
	&=
	(-1)^k
	\sum_{i=0}^{k}
	\binom{m+1}{k-i}\binom{n+1}{i}
	(n-m+k-2i)X_{m+n-k}\\
	&=
	(-1)^k(n-m)\binom{m+n+1}{k}X_{m+n-k}\\
	&=
	\D^{(k)}[L_m,X_n].
\end{align*}
We now consider $a=L_m$ and $b=H_n$.
By \eqref{eq:B-mod-acts-G}, \eqref{eq:Hm-Hn-brac},
and \eqref{eq:Vand-iden}, we have
\begin{align*}
	\sum_{i=0}^{k}
	[\D^{(k-i)}L_m,\D^{(i)}H_n]
	&=
	(-1)^k\sum_{i=0}^{k}
	\binom{m+1}{k-i}\binom{n}{i}(n-i)H_{m+n-k}\\
	&\?+
	(-1)^k\delta_{m+n-k,0}
	\sum_{i=0}^{k}
	\binom{m+1}{k-i}\binom{n}{i}
	(m-k+i)(m-k+i+1)\bc_2\\
	&=
	(-1)^kn\sum_{i=0}^{k}
	\binom{m+1}{k-i}\binom{n-1}{i}H_{m+n-k}\\
	&\?+
	(-1)^k\delta_{m+n-k,0}
	\sum_{i=0}^{k}
	\binom{m+1}{k-i}\binom{n}{i}
(n-i)(n-i-1)\bc_2\\
	&=
	(-1)^kn\binom{m+n}{k}H_{m+n-k}
	+(-1)^k
	n(n-1)\binom{m+n-1}{k}
	\delta_{m+n-k,0}\bc_2\\
	&=\D^{(k)}[L_m,H_n],
\end{align*}
	where the central term vanishes because
	$\delta_{m+n-k,0}\binom{m+n-1}{k}
	=\delta_{m+n-k,0}\binom{k-1}{k}=0$ for $k\in \BZ_+$.
For the remaining nontrivial pairs
$a=H_m$ and $b=H_n,I_n,J_n$,
	\eqref{eq:Lie-B-com} follows directly from
	\eqref{eq:Hm-Hn-brac}--\eqref{eq:H-I;H-J}
	and \eqref{eq:Vand-iden}.
		The cases involving two elements from
		$\{I_m,J_m\mid m\in\BZ\}$ or a central element are immediate.
	Thus, $\mathcal{G}$ is a $\B$-module Lie algebra.
	Consequently, by Lemma \ref{lem:H-mod-Lie-H-mod-alg},
	$U(\mathcal{G})$ is a $\B$-module algebra.
\end{proof}

Set the generating functions
\begin{equation*}
	X(x):=\sum_{n\in \BZ} X_n x^{-n-2},\quad H(x):=\sum_{n\in \BZ}H_nx^{-n-1},
\end{equation*}
where $X=L,I,J$.
Then \eqref{eq:Lm-Ln-brac}--\eqref{eq:H-I;H-J} amount to
\begin{align}
	[L(x_1),L(x_2)]
	&=
	-(\partial_{x_2}^{(1)}L(x_2))
	x_2^{-1}\delta\biggl(\frac{x_1}{x_2}\biggr)
	+2L(x_2)\partial_{x_1}^{(1)}
	\biggl(x_2^{-1}\delta\biggl(\frac{x_1}{x_2}\biggr)\biggr)
	\nonumber\\
	&\?-
	\frac{1}{2}\partial_{x_1}^{(3)}
	\biggl(x_2^{-1}\delta\biggl(\frac{x_1}{x_2}\biggr)\biggr)\bc_1,
	\label{eq:3.12}
	\\
	[L(x_1),H(x_2)]
	&=
	-(\partial_{x_2}^{(1)}H(x_2))
	x_2^{-1}\delta\biggl(\frac{x_1}{x_2}\biggr)
	+H(x_2)\partial_{x_1}^{(1)}
	\biggl(x_2^{-1}\delta\biggl(\frac{x_1}{x_2}\biggr)\biggr)
	\nonumber\\
	&\?+
	2\partial_{x_1}^{(2)}
	\biggl(x_2^{-1}\delta\biggl(\frac{x_1}{x_2}\biggr)\biggr)\bc_2,
	\\
	[L(x_1),I(x_2)]
	&=
	-(\partial_{x_2}^{(1)}I(x_2))
	x_2^{-1}\delta\biggl(\frac{x_1}{x_2}\biggr)
	+2I(x_2)\partial_{x_1}^{(1)}
	\biggl(x_2^{-1}\delta\biggl(\frac{x_1}{x_2}\biggr)\biggr),
	\\
	[L(x_1),J(x_2)]
	&=
	-(\partial_{x_2}^{(1)}J(x_2))
	x_2^{-1}\delta\biggl(\frac{x_1}{x_2}\biggr)
	+2J(x_2)\partial_{x_1}^{(1)}
	\biggl(x_2^{-1}\delta\biggl(\frac{x_1}{x_2}\biggr)\biggr),
	\\
	[H(x_1),H(x_2)]
	&=
	-\partial_{x_1}^{(1)}
	\biggl(x_2^{-1}\delta\biggl(\frac{x_1}{x_2}\biggr)\biggr)\bc_3,
	\\
	[H(x_1),I(x_2)]
	&=
	I(x_2)x_2^{-1}\delta\biggl(\frac{x_1}{x_2}\biggr),
	\;\;
	[H(x_1),J(x_2)]
	=
	-J(x_2)x_2^{-1}\delta\biggl(\frac{x_1}{x_2}\biggr),
	\\
	[I(x_1),I(x_2)]
	&=
	[I(x_1),J(x_2)]
	=
	[J(x_1),J(x_2)]
	=0.
	\label{eq:3.18}
\end{align}
It follows from \cite[(2.7)]{LM} that
\begin{equation}\label{eq:weak-comm}
	(x_1-x_2)^4[a(x_1),b(x_2)]=0
	\quad\text{for }
	a(x),b(x)\in\{L(x),I(x),J(x),H(x)\}.
\end{equation}

Set
\begin{equation*}
	\mathfrak{p}:=
	\bigoplus_{n\ge -1}
	(\BF L_n\oplus\BF I_n\oplus\BF J_n)
	\oplus
	\bigoplus_{n\ge 0}\BF H_n.
\end{equation*}
It follows from \eqref{eq:Lm-Ln-brac}--\eqref{eq:H-I;H-J} that $\mathfrak{p}$ is a
subalgebra of $\mathcal{G}$.
Let $c=(c_1,c_2,c_3)\in \BF^3$, and
let $\BF_c$ be the one-dimensional $\big(\mathfrak{p}\oplus\bigoplus_{i=1}^{3}\BF\bc_i\big)$-module on which $\mathfrak p$ acts trivially and $\bc_i$ acts as the scalar $c_i$ for $1\le i\le3$.
Then we form the induced module
\begin{equation}
	V_{\mathcal{G}}(c,0):=
	U(\mathcal{G})
	\otimes_{U(\mathfrak{p}\oplus\bigoplus_{i=1}^{3}\BF\bc_i)}
	\BF_{c},
\end{equation}
which is a $\mathcal{G}$-module of central charge $c$.

Set $\1:=1\otimes 1$, $v_H:=H(-1)\1$, and $
v_X:=X(-2)\1$ for $X=L,I,J$.
We then obtain the following result.
\begin{theorem}\label{prop:Galilean-VA}
	For any $c=(c_1,c_2,c_3)\in\BF^3$, there exists a
	vertex algebra structure on $V_{\mathcal{G}}(c,0)$,
	uniquely determined by the vacuum vector $\1$
	and the vertex operators
	\begin{equation}\label{eq:ver-ope-2}
		Y(v_X,x)=X(x)=\sum_{n\in\BZ}X(n)x^{-n-2},\;\;
		Y(v_H,x)=H(x)=\sum_{n\in\BZ}H(n)x^{-n-1}
	\end{equation}
	for $X=L,I,J$.
	Moreover, $V_{\mathcal{G}}(c,0)$ is generated by
	$v_L,v_I,v_J,v_H$ as a vertex algebra.
\end{theorem}

\begin{proof}
By \eqref{eq:B-mod-acts-G}, the subspace $\mathfrak{p}+\BF(\bc_1-c_1)+\BF(\bc_2-c_2)+\BF(\bc_3-c_3)$ is a $\B$-submodule of $U(\mathcal{G})$.
Since $U(\mathcal{G})$ is a $\B$-module algebra by Lemma \ref{lem:U(G)-B-mod-alg}, 
the left ideal
\begin{equation}\label{eq:left-ideal}
J_{c}:=	U(\mathcal{G})\big(\mathfrak{p}+\BF(\bc_1-c_1)+\BF(\bc_2-c_2)+\BF(\bc_3-c_3)\big)
\end{equation}
of $U(\mathcal{G})$	is a $\B$-submodule of $U(\mathcal{G})$.	
Hence $U(\mathcal{G})\big/J_{c}$ is naturally a $(U(\mathcal{G}),\B)$-module.
Since $V_{\mathcal{G}}(c,0)\cong U(\mathcal{G})/J_c$ as $U(\mathcal{G})$-modules, 
it follows that $V_{\mathcal{G}}(c,0)$ inherits a $(U(\mathcal{G}),\B)$-module structure, with $\D^{(k)}\1=0$ for $k\in \BZ_+$.
	
Let $T=\{v_L,v_I,v_J,v_H\}$.
Since $V_{\mathcal{G}}(c,0)$ is clearly a locally truncated $\mathcal{G}$-module,
we can define $Y_0(\cdot,x)$ on $T$ by 
\begin{equation*}
	Y_{0}(v_X,x)=X(x),\quad
	Y_{0}(v_H,x)=H(x)\quad \text{for }X=L,I,J.
\end{equation*}
The construction of $V_{\mathcal{G}}(c,0)$ gives $Y_{0}(a,x)\1\in V_{\mathcal{G}}(c,0)[[x]]$ and $\lim\limits_{x\to 0}Y_{0}(a,x)\mathbf{1}=a$ for $a\in T$.
The PBW theorem shows that 
\begin{equation}\label{eq:pbw-basis-span}
	V_{\mathcal{G}}(c,0)=
	\spanf\{
	a^{(1)}(n_1)\cdots a^{(k)}(n_k)\1
	\mid
	a^{(i)}\in\{L,I,J,H\}, k\in \BN, n_i\in \BZ
	\}.
\end{equation}
Comparing coefficients of $x^{-n-1}$ in the defining series for $Y_0$,
we obtain
\begin{equation}\label{eq:coeff-corr}
	(v_X)_n=X(n-1),\quad (v_H)_n=H(n)
\quad\text{for }X=L,I,J,\ n\in\BZ.
\end{equation}
Thus,  
Theorem~\ref{th:LM-general-VA}(iii) follows from \eqref{eq:pbw-basis-span}.
Furthermore, by Lemmas \ref{lem:U(G)-B-mod-alg} and \ref{lem:H-mod-Lie-H-mod-alg}, \eqref{eq:B-mod-acts-G}, and \eqref{eq:taylor}, 
we have
\begin{equation*}
	e^{z\D}Y_0(u,x)e^{-z\D}
	=
	e^{z\D}(Y_{0}(u,x))=Y_{0}(u,x+z)=e^{z\partial_x}Y_0(u,x)\quad\text{for }u\in T
\end{equation*}
on $V_{\mathcal{G}}(c,0)$.
Together with \eqref{eq:weak-comm}, the desired result follows from
Theorem \ref{th:LM-general-VA}.	
	
Moreover, \eqref{eq:pbw-basis-span} and
\eqref{eq:coeff-corr} show that $T$ generates
$V_{\mathcal{G}}(c,0)$ as a vertex algebra.
\end{proof}

As in the case of universal affine and Virasoro vertex algebras (see \cite[Proposition~3.4]{LM} and \cite[Proposition~3.4]{JLM}), 
we have the following module correspondence.

\begin{proposition}\label{prop:mod-equiv}
Every $V_{\mathcal{G}}(c,0)$-module $W$ is a locally truncated $\mathcal{G}$-module of central charge $c=(c_1,c_2,c_3)$, 
with $X(x)=Y_W(v_X,x)$ for $X=L,I,J,H$.
On the other hand, on every locally truncated
$\mathcal{G}$-module $W$ of central charge $c=(c_1,c_2,c_3)$,
there exists a $V_{\mathcal{G}}(c,0)$-module structure $Y_W(\cdot,x)$, 
which is uniquely determined by
\begin{equation*}
	Y_W(v_X,x)=X(x)=\sum_{n\in\BZ}X(n)x^{-n-2},\;\;
	Y_W(v_H,x)=H(x)
	=\sum_{n\in\BZ}H(n)x^{-n-1}
\end{equation*}
for $X=L,I,J$.
\end{proposition}

\begin{proof}
Let $W$ be a $V_{\mathcal{G}}(c,0)$-module.
By \eqref{eq:3.12}--\eqref{eq:3.18} and \cite[Lemma~2.16]{LM},
the defining commutator relations of $\mathcal{G}$ hold on $W$,
with $\bc_i$ acting as $c_i\operatorname{id}_W$ for $i=1,2,3$.
Thus, $W$ is a $\mathcal{G}$-module of central charge $c=(c_1,c_2,c_3)$.
Since $X(x)W=Y_{W}(v_{X},x)W\subset W((x))$ for $X=L,I,J,H$, $W$ is a locally truncated $\mathcal{G}$-module.

Conversely, let $W$ be a locally truncated $\mathcal{G}$-module of central charge $c=(c_1,c_2,c_3)$.
	Set
	$S:=\{L(x),I(x),J(x),H(x)\}\subset\mathcal{E}(W)$.
	By \eqref{eq:weak-comm} and Theorem \ref{th:local-set-VA},
	$\langle S\rangle$ is a vertex algebra and $W$ is a faithful
	$\langle S\rangle$-module, with $Y_W(X(x),z)=X(z)$ for $X(x)\in S$.
	By \eqref{eq:3.12}--\eqref{eq:3.18} and \cite[Lemma 2.16]{LM},
	$\langle S\rangle$ is a locally truncated $\mathcal{G}$-module of central charge $c$,
	where $L(n)$, $I(n)$, $J(n)$, and $H(n)$ are given by
	$L(x)_{n+1}$, $I(x)_{n+1}$, $J(x)_{n+1}$, and $H(x)_n$,
	respectively.
	By the creation property (see Definition \ref{def:VA}(3)),
	\begin{equation*}
		X(n)\operatorname{id}_W=0
		\quad\text{for }X=L,I,J\text{ and }n\ge-1,
		\quad
		\text{and }H(m)\operatorname{id}_W=0
		\quad\text{for }m\ge0.
	\end{equation*}
Thus, by the universal property of $V_{\mathcal{G}}(c,0)$,
there exists a unique $\mathcal{G}$-module homomorphism
$\psi:V_{\mathcal{G}}(c,0)\to\langle S\rangle$
such that $\psi(\1)=\operatorname{id}_W$.	
Since $\psi$ is a $\mathcal{G}$-module homomorphism, we have
\begin{equation*}
	\psi((v_{X^1})_{n_1}\cdots (v_{X^k})_{n_k}\1)=X^1(x)_{n_1}\cdots X^k(x)_{n_k}\operatorname{id}_{W}
\end{equation*}
for
$v_{X^i}\in\{v_L,v_I,v_J,v_H\}$, $n_i\in\BZ$, and $k\in \BN$.
Therefore, by Theorem \ref{th:gen-con-V-mod} with
$T=\{v_L,v_I,v_J,v_H\}$,
$W$ is a $V_{\mathcal{G}}(c,0)$-module with $Y_W(v,x)=\psi(v)$
for $v\in V_{\mathcal{G}}(c,0)$.
\end{proof}

Recall that (see, for example, \cite[Definition~2.6]{JLM}) an {\em $\BN$-graded vertex algebra} is a vertex algebra $V$ equipped with an
$\BN$-grading $V=\bigoplus\limits_{n\in\BN}V_{(n)}$ such that
$\1\in V_{(0)}$ and
\begin{equation*}
	u_kV_{(n)}\subset V_{(m+n-k-1)}
	\quad\text{for }u\in V_{(m)},\ m,n\in\BN,\text{ and }k\in\BZ,
\end{equation*}
where $V_{(n)}=0$ for $n<0$.
	For an $\BN$-graded vertex algebra $V$,
	an {\em $\BN$-graded $V$-module} is a $V$-module
	$W=\bigoplus\limits_{n\in\BN}W_{(n)}$ such that
\begin{equation*}
		u_kW_{(n)}\subset W_{(m+n-k-1)}
		\quad\text{for }u\in V_{(m)}, m,n\in\BN,\text{ and } k\in\BZ,
\end{equation*}
	where $W_{(n)}=0$ for $n<0$.
	Following \cite{JLM},
	two $\BN$-graded $V$-modules $W^1$ and $W^2$
	are {\em equivalent} if there exist a $V$-module isomorphism
	$f:W^1\to W^2$ and $s\in\BZ$ such that
\begin{equation*}
		f(W^1_{(n)})=W^2_{(n+s)}
		\quad\text{for }n\in\BZ.
\end{equation*}
	Every nonzero $\BN$-graded $V$-module is therefore
	equivalent to one with nonzero degree-zero component.

%

We make $\mathcal{G}$ a $\BZ$-graded Lie algebra by
setting
\begin{equation}\label{eq:deg-bas-G}
	\deg L_n=\deg I_n=\deg J_n=\deg H_n=-n,\quad
	\deg\bc_i=0
\end{equation}
for $n\in\BZ$ and $i=1,2,3$.
Since $J_c$ is a graded left ideal of $U(\mathcal{G})$,
the PBW theorem gives an $\BN$-grading on
$V_{\mathcal{G}}(c,0)(\cong U(\mathcal{G})/J_c)$:
\begin{equation}\label{eq:Galilean-VOA-grading}
	V_{\mathcal{G}}(c,0)
	=\bigoplus_{n\in\BN}V_{\mathcal{G}}(c,0)_{(n)},
\end{equation}
where $V_{\mathcal{G}}(c,0)_{(n)}$ has a basis consisting of
\begin{equation}\label{eq:Galilean-VOA-PBW}
	L(-m_1)\cdots L(-m_r)
	I(-n_1)\cdots I(-n_s)
	J(-k_1)\cdots J(-k_t)
	H(-l_1)\cdots H(-l_q)\1
\end{equation}
for $r,s,t,q\in\BN$,
$m_1\ge\ldots\ge m_r\ge2$,
$n_1\ge\ldots\ge n_s\ge2$,
$k_1\ge\ldots\ge k_t\ge2$,
$l_1\ge\ldots\ge l_q\ge1$, 
and $\sum_{i=1}^{r}m_i
+\sum_{i=1}^{s}n_i
+\sum_{i=1}^{t}k_i
+\sum_{i=1}^{q}l_i=n$.
In particular,
$V_{\mathcal{G}}(c,0)_{(0)}=\BF\1$ and
$\dim V_{\mathcal{G}}(c,0)_{(n)}<\infty$ for $n\in\BN$.
Set $V_{\mathcal{G}}(c,0)_{(n)}=0$ for $n<0$.

By \eqref{eq:coeff-corr},
\begin{equation*}
	(v_X)_rV_{\mathcal{G}}(c,0)_{(n)}
	\subset
	V_{\mathcal{G}}(c,0)_{(\deg v_X+n-r-1)}
\end{equation*}
for $X=L,I,J,H$, $r\in\BZ$, and $n\in\BN$.
Since $v_L,v_I,v_J,v_H$ generate $V_{\mathcal{G}}(c,0)$ by Theorem \ref{prop:Galilean-VA},
it follows from the iterate formula (see \cite[(3.1.12)]{LL}) that
$V_{\mathcal{G}}(c,0)$ is an $\BN$-graded vertex algebra.

We use the definition of a vertex operator algebra
over $\BF$ given by Dong and Ren in
\cite[Definition~2.3]{DR} (cf. \cite{FHL,LL}).
A direct calculation gives the following result.

\begin{proposition}\label{prop:Galilean-VOA}
	For any $c=(c_1,c_2,c_3)\in\BF^3$,
	the vertex algebra $V_{\mathcal{G}}(c,0)$,
	equipped with the grading \eqref{eq:Galilean-VOA-grading}
	and the conformal vector $\omega=-L(-2)\1$,
	is a vertex operator algebra of central charge $c_1$.
\end{proposition}

Next, we study quotient vertex algebras of $V_{\mathcal{G}}(c,0)$.
Recall that a {\em restricted Lie algebra} over $\BF$ is a Lie algebra
$\fg$ equipped with a mapping $a\mapsto a^{[p]}$, called the
{\em $p$-mapping}, satisfying the following conditions:
\begin{enumerate}[(i)]
	\item $\operatorname{ad} a^{[p]}=(\operatorname{ad}a)^{p}$ for $a\in \fg$.
	\item $(\lambda a)^{[p]}=\lambda^{p}a^{[p]}$ for $\lambda\in \BF$, $a\in\fg$.
	\item $(a+b)^{[p]}=a^{[p]}+b^{[p]}+\sum_{i=1}^{p-1}s_{i}(a,b)$ for $a,b\in \fg$, 
	where $s_{i}(a,b)$ are given by 
	\begin{equation*}
		(\operatorname{ad}(a\otimes X+b\otimes 1))^{p-1}(a\otimes 1)
		=
		\sum_{i=1}^{p-1}is_i(a,b)\otimes X^{i-1}
	\end{equation*}
	in the Lie algebra $\fg\otimes \BF[X]$ over the polynomial ring $\BF[X]$.
\end{enumerate}

The following result is well known in modular Lie algebra theory.
\begin{theorem}[{\cite{Jacobson1962}; or \cite[Chapter~2, Theorem~2.3]{SF}}]\label{th:4.8}
	Let $\fg$ be a Lie algebra with a basis $\{a_i\mid i\in I\}$ for some
	index set $I$ such that there are $b_i\in \fg$ with
	$(\ad a_i)^{p}=\ad b_i$ for $i\in I$.
	Then there exists a restricted Lie algebra structure on $\fg$ with
	the $p$-mapping uniquely determined by $a_i^{[p]}=b_i$ for $i\in I$.
\end{theorem}

For any restricted Lie algebra $\fg$, it is well known that
$a^{p}-a^{[p]}$ lies in the center of the universal enveloping algebra
$U(\fg)$ for every $a\in\fg$.
The {\em restricted universal enveloping algebra}
$\mathfrak{u}(\fg)$ is defined to be the quotient of $U(\fg)$
by the ideal generated by the elements
$a^{p}-a^{[p]}$ for $a\in\fg$.

Following \cite{JLM}, for $m\in \BZ$, set 
\begin{equation*}
	\delta_{p\mid m}:=
	\begin{cases}
		1 & \text{if }p\mid m,\\
		0 & \text{if }p\nmid m.
	\end{cases}
\end{equation*}
Recall that Fermat's little theorem gives $	n^{p}\equiv n \pmod{p}$ for $n\in\BZ$.
Consequently, we have 
\begin{equation}\label{eq:app-fermat-x}
	\prod_{j=0}^{p-1}(x-j)=x^p-x\quad \text{in } \BF[x].
\end{equation}
In particular, 
for $m,n\in\BZ$,
\begin{equation}\label{eq:app-fermat}
	\prod_{j=0}^{p-1}(n+jm)
	=\delta_{p\mid m}n
	\quad\text{in }\BF.
\end{equation}

\begin{lemma}\label{lem:rest-stru-on-G}
	There is a restricted Lie algebra structure on
	$\mathcal{G}$, whose $p$-mapping is uniquely determined by
	\begin{equation}\label{eq:w22-p-map}
		\bc_i^{[p]}=\bc_i,
		\quad
		L_m^{[p]}=\delta_{p\mid m}L_{pm},
		\quad
		I_m^{[p]}=
		J_m^{[p]}=0,\quad
		H_m^{[p]}=H_{pm}
	\end{equation}
	for $i=1,2,3$ and $m\in\BZ$.
\end{lemma}

\begin{proof}
Let $m,n\in\BZ$.
By \cite[Lemma~3.8]{JLM}, applied to the Virasoro subalgebra with generators $-L_n$, we have
\begin{equation}\label{eq:L}
	(\ad L_m)^pL_n
	=\delta_{p\mid m}[L_{pm},L_n].
\end{equation}
Using \eqref{eq:app-fermat} and \eqref{eq:Hm-Hn-brac}--\eqref{eq:L-I;L-J},
we obtain
\begin{align}
	(\ad L_m)^pX_n
	&=
	\Big(\prod_{j=0}^{p-1}(n+(j-1)m)\Big)X_{pm+n}\nonumber\\
	&=
	\delta_{p\mid m}nX_{pm+n}
	=
	\delta_{p\mid m}[L_{pm},X_n]
	\quad\text{for }X=I,J,\label{eq:IJ}
	\\
	(\ad L_m)^pH_n
	&=
	\Big(\prod_{j=0}^{p-1}(n+jm)\Big)H_{pm+n}+
	m(m+1)\delta_{pm+n,0}
	\Big(\prod_{j=0}^{p-2}(n+jm)\Big)\bc_2\nonumber\\
	&=
	\delta_{p\mid m}nH_{pm+n}
	=\delta_{p\mid m}[L_{pm},H_n],\label{eq:H}
\end{align}
	where the central term vanishes because
	$\delta_{pm+n,0}\ne0$ implies $n=0$ in $\BF$.
Since $\bc_i$ are central, \eqref{eq:L}--\eqref{eq:H} imply
	\begin{equation}\label{eq:p-map-Lm}
		(\ad L_m)^p
		=
		\ad(\delta_{p\mid m}L_{pm})
		\quad\text{on }\mathcal{G}.
	\end{equation}
	
For the remaining basis elements, a direct calculation from
\eqref{eq:Hm-Hn-brac}--\eqref{eq:H-I;H-J}, using $p>2$, gives
\begin{equation}\label{eq:p-map-other}
	(\ad H_m)^p=\ad H_{pm},\qquad
	(\ad I_m)^p=(\ad J_m)^p=0,\qquad
	(\ad\bc_i)^p=\ad\bc_i=0
\end{equation}
for $i=1,2,3$.
The result follows from Theorem \ref{th:4.8},
\eqref{eq:p-map-Lm}, and \eqref{eq:p-map-other}.
\end{proof}

Lemma \ref{lem:rest-stru-on-G} gives the following corollary.
\begin{corollary}\label{coro:p-cent}
	For $n\in\BZ$, the elements
	$L_n^p-\delta_{p\mid n}L_{pn}$, $I_n^p$, $J_n^p$,
	and $H_n^p-H_{pn}$ lie in the center of $U(\mathcal{G})$.
\end{corollary}

For later use, we compute the following vertex operators on a $V_{\mathcal{G}}(c,0)$-module.

\begin{proposition}\label{prop:vertex-for-generator}
	Let $W$ be a $V_{\mathcal{G}}(c,0)$-module.
	The following relations hold on $W$ for $n\ge 2$, $m\ge 1$,
	and $X=I,J$:
	\begin{align}
		&\?Y_W\big((L(-n)^{p}-\delta_{p\mid n}L(-pn))\1,x\big)\nonumber\\
		&=\sum_{j\in \BN}(-1)^{j}\binom{-n+1}{j}\big(L(-n-j)^{p}-\delta_{p\mid (n+j)}L(-p(n+j))\big)x^{jp}+\sum_{j\in \BN}(-1)^{n-j}\binom{-n+1}{j}\nonumber\\
		&\?\times \big(L(j-1)^{p}-\delta_{p\mid (j-1)}L((j-1)p)\big)x^{(-n+1-j)p},\label{eq:Y-L-n-p}\\
		&\?Y_W(X(-n)^{p}\1,x)\nonumber\\
	&=\sum_{j\in\BN}(-1)^j \binom{-n+1}{j}X(-n-j)^p x^{jp}
	+\sum_{j\in\BN}(-1)^{n-j}\binom{-n+1}{j}
	X(j-1)^p x^{(-n+1-j)p},\label{eq:Y-IJ-n-p}\\
		&\?Y_W\big((H(-m)^{p}-H(-pm))\1,x\big)\nonumber\\
		&=\sum_{j\in\BN}\binom{m+j-1}{j}\big(H(-m-j)^{p}-H(-p(m+j))\big)x^{jp}+(-1)^{m-1}\sum_{j\in\BN}\binom{m+j-1}{j}\nonumber\\
		&\?\times\big(H(j)^{p}-H(pj)\big)x^{-p(m+j)}.\label{eq:Y-H-mp}
	\end{align}
\end{proposition}

\begin{proof}
Let $n\ge2$ and $m\ge1$.
With $-L_k$ for $k\in\BZ$ as the Virasoro generators,
\eqref{eq:Y-L-n-p} follows from the calculation
in the proof of \cite[Proposition~3.12]{JLM}.
	
For $X=I,J$, we have $(v_X)_k=X(k-1)$ for $k\in\BZ$.
Since $[X_r,X_s]=0$ for $r,s\in\BZ$, 
the operators $(v_X)_k$ for $k\in\BZ$ commute pairwise.
Thus, \cite[Lemma~5.6]{LM}
gives \eqref{eq:Y-IJ-n-p}.

Since $(v_H)_k=H(k)$ for $k\in\BZ$, the operators
$(v_H)_k$ for $k\in\BN$ mutually commute, and so do
$(v_H)_{-k}$ for $k\in\BZ_+$.
Hence, by \cite[Lemma~5.6]{LM},
\begin{align}\label{eq:H-m-p}
	&\?Y_W(H(-m)^p\1,x)\nonumber\\
	&=
	\sum_{j\in\BN}\binom{m+j-1}{j}
	H(-m-j)^p x^{jp}
	+
	(-1)^{m-1}
	\sum_{j\in\BN}\binom{m+j-1}{j}
	H(j)^p x^{-p(m+j)}.
\end{align}
Furthermore,
by the iterate formula (see \cite[(3.1.11)]{LL}), we obtain
\begin{align}\label{eq:H-pm}
	Y_W(H(-pm)\1,x)
	&=
	\sum_{j\in\BN}\binom{-pm}{j}(-1)^j
	\big(
	H(-pm-j)x^j
	-
	(-1)^{pm}H(j)x^{-pm-j}
	\big)
	\nonumber\\
	&=
	\sum_{j\in\BN}\binom{pm+j-1}{j}
	\big(
	H(-pm-j)x^j
	+
	(-1)^{m-1}H(j)x^{-pm-j}
	\big)\nonumber\\
	&=
	\sum_{j\in\BN}\binom{m+j-1}{j}
	\big(
	H(-p(m+j))x^{jp}
	+
	(-1)^{m-1}
	H(pj)x^{-p(m+j)}
	\big),
\end{align}
where the last equality follows from \cite[Lemma~5.4]{JLM}.
	Combining \eqref{eq:H-pm} with \eqref{eq:H-m-p}
	gives \eqref{eq:Y-H-mp}.
\end{proof}

By Lemma \ref{lem:rest-stru-on-G}, the subspace
\begin{equation*}
	\fg_{\mathcal{G}}
	:=\BF L_0\oplus\BF I_0\oplus\BF J_0\oplus\BF H_0
\end{equation*}
is a graded restricted Lie subalgebra of $\mathcal{G}$
with nontrivial brackets
\begin{equation}\label{eq:non-bra-g}
	[H_0,I_0]=I_0,\quad [H_0,J_0]=-J_0
\end{equation}
and $p$-mapping determined by
\begin{equation}\label{eq:rest-stru-g}
	L_0^{[p]}=L_0,\quad
	I_0^{[p]}=J_0^{[p]}=0,\quad
	H_0^{[p]}=H_0.
\end{equation}
Let $\chi$ be a linear functional on $\fg_{\mathcal{G}}$,
and let $I^\chi$ be the $\mathcal{G}$-submodule of
$V_{\mathcal{G}}(c,0)$ generated by
\begin{align}\label{eq:gen-I-chi}
	&\big(L(-n)^p-\delta_{p\mid n}L(-pn)
	-\delta_{n,2}\chi(L_0)^p\big)\1,
	\nonumber\\
	&\big(X(-n)^p-\delta_{n,2}\chi(X_0)^p\big)\1,
	\nonumber\\
	&\big(H(-m)^p-H(-pm)
	-\delta_{m,1}\chi(H_0)^p\big)\1
\end{align}
for $X=I,J$, $n\ge2$, and $m\ge1$.

\begin{theorem}\label{th:I-mu-iedal}
The $\mathcal{G}$-submodule $I^\chi$ is an ideal of $V_{\mathcal{G}}(c,0)$, 
which equals the ideal generated by the four vectors
	\begin{equation}\label{eq:ideal-gener}
		\bigl(L(-2)^p-\chi(L_0)^p\bigr)\1,\quad
		\bigl(X(-2)^p-\chi(X_0)^p\bigr)\1,\quad
		\bigl(H(-1)^p-H(-p)-\chi(H_0)^p\bigr)\1
	\end{equation}
for $X=I,J$.
\end{theorem}

\begin{proof}
Since $I^\chi$ is a $\mathcal{G}$-submodule of $V_{\mathcal{G}}(c,0)$,
it is a left ideal of $V_{\mathcal{G}}(c,0)$.
By \cite[Remark~3.9.8]{LL} (see also \cite[Section~2]{JLM}),
it is sufficient to show that $I^{\chi}$ is $\B$-stable.
	
By Proposition \ref{prop:vertex-for-generator} and \eqref{eq:B-mod-str},
for $n\ge2$, $m\ge1$, and $k\in\BN$,
we obtain
\begin{align}
	&\?\D^{(k)}
	\big((L(-n)^p-\delta_{p\mid n}L(-pn)
	-\delta_{n,2}\chi(L_0)^p)\1\big)
	\nonumber\\
	&=
	\begin{cases}
		(-1)^r\binom{-n+1}{r}
		\big(
		L(-n-r)^p-\delta_{p\mid(n+r)}L(-p(n+r))
		-\delta_{n+r,2}\chi(L_0)^p
		\big)\1,
		&\text{if }k=pr,\\
		0,&\text{otherwise},
	\end{cases}
	\label{eq:D-L-p-center}\\
	&\?\D^{(k)}
	\big((X(-n)^p-\delta_{n,2}\chi(X_0)^p)\1\big)
	\nonumber\\
	&=
	\begin{cases}
		(-1)^r\binom{-n+1}{r}
		\big(X(-n-r)^p
		-\delta_{n+r,2}\chi(X_0)^p\big)\1,
		&\text{if }k=pr,\\
		0,&\text{otherwise},
	\end{cases}
	\label{eq:D-IJ-p-center}\\
	&\?\D^{(k)}
	\big((H(-m)^p-H(-pm)
	-\delta_{m,1}\chi(H_0)^p)\1\big)
	\nonumber\\
	&=
	\begin{cases}
		\binom{m+r-1}{r}
		\big(H(-m-r)^p-H(-p(m+r))
		-\delta_{m+r,1}\chi(H_0)^p\big)\1,
		&\text{if }k=pr,\\
		0,&\text{otherwise},
	\end{cases}
	\label{eq:D-H-p-center}
\end{align}
where $X=I,J$.
Since $V_{\mathcal{G}}(c,0)$ is a $(U(\mathcal{G}),\B)$-module,
it follows from \eqref{eq:A-H} that $I^\chi$ is a $\B$-submodule of $V_{\mathcal{G}}(c,0)$.
Therefore, $I^{\chi}$ is an ideal of $V_{\mathcal{G}}(c,0)$.
	
Taking $n=2$ and $m=1$ in \eqref{eq:D-L-p-center}--\eqref{eq:D-H-p-center},
we obtain all the generators in \eqref{eq:gen-I-chi} from the four vectors in \eqref{eq:ideal-gener} by applying $\D^{(pr)}$ for $r\in\BN$.
Since $I^\chi$ is an ideal, 
these four vectors  in \eqref{eq:ideal-gener} generate $I^\chi$ as an ideal.
\end{proof}

For $\chi\in\fg_{\mathcal{G}}^*$, set
\begin{equation}
	V^\chi_{\mathcal{G}}(c,0)
	:=V_{\mathcal{G}}(c,0)/I^\chi,
\end{equation}
which is naturally a $\mathcal{G}$-module and a vertex algebra.
We have the following characterization of its modules.

\begin{proposition}\label{prop:quot-mod-equiv}
	Every $V^\chi_{\mathcal{G}}(c,0)$-module $W$ is a locally
	truncated $\mathcal{G}$-module of central charge $c$ such that
	\begin{align}
		L(-j-2)^p-\delta_{p\mid(j+2)}L(-p(j+2))
		&=\delta_{j,0}\chi(L_0)^p\operatorname{id}_W,\label{eq:L-2}\\
		X(-j-2)^p
		&=\delta_{j,0}\chi(X_0)^p\operatorname{id}_W,\label{eq:IJ-2}\\
		H(-j-1)^p-H(-p(j+1))
		&=\delta_{j,0}\chi(H_0)^p\operatorname{id}_W\label{eq:H-1}
	\end{align}
	on $W$ for $X=I,J$ and $j\in\BZ$.
	Furthermore, the converse is also true.
\end{proposition}

\begin{proof}
Note that a $V^{\chi}_{\mathcal{G}}(c,0)$-module is naturally a $V_{\mathcal{G}}(c,0)$-module.
By Proposition \ref{prop:mod-equiv}, a $V_{\mathcal{G}}(c,0)$-module exactly amounts to a locally truncated $\mathcal{G}$-module of central charge $c$.
Let $W$ be a $V_{\mathcal{G}}(c,0)$-module.
Applying Proposition~\ref{prop:vertex-for-generator}
with $n=2$ and $m=1$ and reindexing the sums with negative powers
of $x$, we obtain
\begin{align}
	Y_{W}\big((L(-2)^p-\chi(L_0)^p)\1,x\big)
	&=
	\sum_{j\in\BZ}
	\Big(
	L(-j-2)^p
	-\delta_{p\mid(j+2)}L(-p(j+2))\nonumber\\
	&\?-\delta_{j,0}\chi(L_0)^p\operatorname{id}_{W}
	\Big)x^{jp},
	\label{eq:ver-ope-L}\\
	Y_{W}\big((X(-2)^p-\chi(X_0)^p)\1,x\big)&=
	\sum_{j\in\BZ}
	\big(
	X(-j-2)^p
	-\delta_{j,0}\chi(X_0)^p\operatorname{id}_{W}
	\big)x^{jp},\label{eq:ver-ope-gen-IJ}\\
	Y_{W}\big((H(-1)^p-H(-p)-\chi(H_0)^p)\1,x\big)
	&=
	\sum_{j\in\BZ}
	\Big(
	H(-j-1)^p
	-H(-p(j+1))\nonumber\\
	&\?-\delta_{j,0}\chi(H_0)^p\operatorname{id}_{W}
	\Big)x^{jp},\label{eq:ver-ope-H-4}
\end{align}
where $X=I,J$.

If $W$ is a $V^{\chi}_{\mathcal{G}}(c,0)$-module, since the elements in \eqref{eq:ideal-gener} are in $I^{\chi}$, the left-hand sides of equations \eqref{eq:ver-ope-L}--\eqref{eq:ver-ope-H-4} vanish.
Comparing coefficients gives \eqref{eq:L-2}--\eqref{eq:H-1} for all $j\in \BZ$.
The converse is also straightforward from Proposition \ref{prop:mod-equiv}, Theorem \ref{th:I-mu-iedal}, and \eqref{eq:ver-ope-L}--\eqref{eq:ver-ope-H-4}.
\end{proof}

By Theorem~\ref{th:I-mu-iedal}, $I^0$ is generated
as an ideal by homogeneous vectors in \eqref{eq:ideal-gener} for $\chi=0$, and hence is graded.
On the other hand,
by the definition of $I^0$ (see \eqref{eq:gen-I-chi} for $\chi=0$) and Corollary \ref{coro:p-cent},
$I^0$ is spanned by
\begin{equation}\label{eq:I0-central-span}
	(L(-n)^p-\delta_{p\mid n}L(-pn))u,\quad
	X(-n)^pu,\quad
	(H(-m)^p-H(-pm))u
\end{equation}
for $X=I,J$, $n\ge2$, $m\ge1$,
and $u\in V_{\mathcal{G}}(c,0)$.
Consequently, $\1 \notin I^0$.
Thus, $I^0$ is a nonzero proper ideal, 
and $V_{\mathcal{G}}(c,0)$ is not simple.
The quotient vertex algebra $V^0_{\mathcal{G}}(c,0)$
is naturally $\BN$-graded with $(V^0_{\mathcal{G}}(c,0))_{(0)}=\BF\1$.

\section{Irreducible $\BN$-graded modules for $V_{\mathcal{G}}(c,0)$ and $V^0_{\mathcal{G}}(c,0)$}\label{sec:4}

In this section, we determine the Zhu algebras of
$V_{\mathcal{G}}(c,0)$ and $V^0_{\mathcal{G}}(c,0)$
and use them to classify the irreducible $\BN$-graded
modules of these vertex algebras.

To begin with, we construct $\BN$-graded
$V_{\mathcal{G}}(c,0)$-modules from $\fg_{\mathcal{G}}$-modules.
Let $\mathcal{G}_{(n)}$ denote the degree-$n$ homogeneous
subspace of $\mathcal{G}$, and set
$\mathcal{G}_-=\bigoplus\limits_{n\le0}\mathcal{G}_{(n)}$.
Fix $c=(c_1,c_2,c_3)\in\BF^3$.
For a $\fg_{\mathcal{G}}$-module $U$,
place $U$ in degree zero and extend the $\fg_{\mathcal{G}}$-action to $\mathcal{G}_-$
by letting $\mathcal{G}_{(n)}$ act trivially for $n<0$,
and $\bc_i$ act as the scalar $c_i$ for $i=1,2,3$.
Consider the induced module
\begin{equation}\label{eq:def-M-U}
	M_{\mathcal{G}}(c,U):=
	U(\mathcal{G})\otimes_{U(\mathcal{G}_-)}U.
\end{equation}
By the PBW theorem, $M_{\mathcal{G}}(c,U)$ admits an
$\BN$-grading
$M_{\mathcal{G}}(c,U)
=\bigoplus\limits_{n\in\BN}M_{\mathcal{G}}(c,U)_{(n)}$
such that $M_{\mathcal{G}}(c,U)_{(0)}=U$  and
\begin{equation*}
 X(m)M_{\mathcal{G}}(c,U)_{(n)}\subset M_{\mathcal{G}}(c,U)_{(n-m)}
\end{equation*}
for $X=L,I,J,H$, $m\in\BZ$, and $n\in\BN$,
where $M_{\mathcal{G}}(c,U)_{(m)}=0$ for $m<0$.
Thus, $M_{\mathcal{G}}(c,U)$ is a locally truncated
$\mathcal{G}$-module of central charge $c$.
By Proposition \ref{prop:mod-equiv},
$M_{\mathcal{G}}(c,U)$ is a $V_{\mathcal{G}}(c,0)$-module.
Since $V_{\mathcal{G}}(c,0)$ is generated by
$v_L,v_I,v_J,v_H$ and
$(v_X)_{\deg v_X-1+m}=X(m)$
for $X=L,I,J,H$ and $m\in\BZ$,
$M_{\mathcal{G}}(c,U)$ is an $\BN$-graded
$V_{\mathcal{G}}(c,0)$-module.

Let $L_{\mathcal{G}}(c,U)$ be the quotient of
	$M_{\mathcal{G}}(c,U)$ by the sum of all graded
	$\mathcal{G}$-submodules with zero degree-zero component.
Then $L_{\mathcal{G}}(c,U)$ is an $\BN$-graded module with
\begin{equation}\label{eq:L-0=U}
L_{\mathcal{G}}(c,U)_{(0)}=U.
\end{equation}
It can be readily seen that if $U$ is an irreducible
$\fg_{\mathcal{G}}$-module, then $L_{\mathcal{G}}(c,U)$ is an
irreducible $\BN$-graded $\mathcal{G}$-module, and also an irreducible
$\BN$-graded $V_{\mathcal{G}}(c,0)$-module.

We next recall the Zhu algebra $A(V)$ associated with an 
$\BN$-graded vertex algebra $V$ (see \cite{DR}, \cite{FZ}, \cite{Zhu1996}).
Following \cite{Zhu1996}, for $u,v\in V$ with $u$ homogeneous, 
define
\begin{equation}\label{eq:u*v}
	u*v=\Res_{x} \frac{(1+x)^{\deg u}}{x} Y(u,x)v,
\end{equation}
and extend linearly in $u$.
Then $(V,*)$ carries a nonassociative algebra structure.
Following \cite{DR}, 
for $u,v\in V$ with $u$ homogeneous and for $k\in \BN$, define
\begin{equation}\label{eq:u-circ-v}
	u\circ_{k}v=\Res_{x} \frac{(1+x)^{\deg u}}{x^{k+2}} Y(u,x)v,
\end{equation}
extend linearly in $u$, and set
\begin{equation}\label{eq:def-O(V)}
	O(V):=\spanf\{u\circ_{k} v\mid u,v\in V, k\in \BN\}.
\end{equation}

\begin{proposition}[\cite{DR}]
	The subspace $O(V)$ is an ideal of the nonassociative algebra $(V,*)$, 
	and the quotient $A(V):=V/O(V)$ is an associative algebra, 
	where $\mathbf{1}+O(V)$ is the identity.
\end{proposition}
We write $[u]=u+O(V)\in A(V)$ for $u\in V$.
The product in $A(V)$ is written as $[u][v]=[u*v]$
for $u,v\in V$.

\begin{theorem}\label{th:A(V)}
The Zhu algebra $A(V_{\mathcal{G}}(c,0))$ is isomorphic to $U(\fg_{\mathcal{G}})$ as algebras.
\end{theorem}

\begin{proof}
Define a linear map $\psi:\fg_{\mathcal{G}}\to A(V_{\mathcal{G}}(c,0))$ by $\psi(X_0)=[v_X]$ for $X=L,I,J,H$.
We first show that $\psi$ is a Lie algebra homomorphism.

By \cite[Lemma~3.2(3)]{DR}, for $u,v\in V_{\mathcal{G}}(c,0)$
with $u$ homogeneous, we have
\begin{equation}\label{eq:zhu-comm}
	[u][v]-[v][u]
	=
	\sum_{j\in\BN}
	\binom{\deg u-1}{j}[u_jv].
\end{equation}
For $X=I,J,H$, using \eqref{eq:zhu-comm},
\eqref{eq:Hm-Hn-brac}--\eqref{eq:L-I;L-J},
and \cite[Lemma~2.9]{JLM}, we obtain
\begin{equation*}
	[v_L][v_X]-[v_X][v_L]
	=[(v_L)_0v_X+(v_L)_1v_X]
	=[-\D^{(1)}v_X-(\deg v_X)v_X]=0.
\end{equation*}
Next, since $\deg v_H=1$,
\eqref{eq:zhu-comm}, \eqref{eq:coeff-corr},
and \eqref{eq:H-I;H-J} give
\begin{align*}
	[v_H][v_I]-[v_I][v_H]
	&=[H(0)v_I]=[v_I],\\
	[v_H][v_J]-[v_J][v_H]
	&=[H(0)v_J]=-[v_J].
	\nonumber
\end{align*}	
Finally, since $(v_I)_jv_J=0$ for $j\ge0$, 
\eqref{eq:zhu-comm} also gives
	\begin{equation*}
		[v_I][v_J]-[v_J][v_I]=0.
	\end{equation*}
These relations agree with \eqref{eq:non-bra-g}, so $\psi$ is a Lie algebra homomorphism.
	By the universal property of $U(\fg_{\mathcal{G}})$,
	$\psi$ extends uniquely to an algebra homomorphism
	\begin{equation*}
		\bar\psi:U(\fg_{\mathcal{G}})
		\to A(V_{\mathcal{G}}(c,0)).
	\end{equation*}
	
We next show that $\bar\psi$ is surjective.
Since $V_{\mathcal{G}}(c,0)$ is spanned by homogeneous elements
\begin{equation}\label{eq:pbw-basis-for V}
	(v_{X^1})_{-n_1}\cdots(v_{X^k})_{-n_k}\1,
\end{equation}
where $k\in\BN$, $n_i\in\BZ_+$, and
$X^i\in\{L,I,J,H\}$ for $1\le i\le k$,
it suffices to show by induction on $N\in\BN$ that $[u]\in\operatorname{Im}\bar\psi$ for every homogeneous $u$ of degree $N$.
The assertion is clear when $N=0$,
since $u$ is a scalar multiple of $\1$.
Let $N>0$, and assume that the assertion holds for all degrees less than $N$.
By linearity, we may take $u$ of the form \eqref{eq:pbw-basis-for V} and write
\begin{equation*}
	u=a_{-n}v,
\end{equation*}
where $a\in\{v_L,v_I,v_J,v_H\}$,
$n\ge1$, and $v$ is homogeneous.	
	
Consider first $n\ge2$.
By \eqref{eq:u-circ-v} and \eqref{eq:def-O(V)},
$a\circ_{n-2}v
=a_{-n}v+\sum\limits_{j\ge1}\binom{\deg a}{j}a_{j-n}v
\in O(V_{\mathcal{G}}(c,0))$.
Hence
\begin{equation}\label{eq:a--nv-in-A(V)}
	[u]=[a_{-n}v]
	=
	-\sum_{j\ge1}
	\binom{\deg a}{j}[a_{j-n}v].
\end{equation}
For $n=1$, \eqref{eq:u*v} gives $a*v=a_{-1}v+\sum\limits_{j\ge1}\binom{\deg a}{j}a_{j-1}v$,
and thus
\begin{equation}\label{eq:reduce-minus-one}
	[u]=[a_{-1}v]
	=
	[a][v]
	-
	\sum_{j\ge1}
	\binom{\deg a}{j}[a_{j-1}v].
\end{equation}
In both cases, $\deg v=N-\deg a-n+1<N$,
	and every nonzero term $a_{j-n}v$ with $j\ge1$
	has degree $N-j<N$.
The induction hypothesis therefore gives
$[u]\in\operatorname{Im}\bar\psi$.
This proves that $\bar\psi$ is surjective.

It remains to prove that $\bar\psi$ is injective.
Take $U=U(\fg_{\mathcal{G}})$, regarded as a
$\fg_{\mathcal{G}}$-module under left multiplication.
By \cite[Proposition~2.10]{JLM},
$M_{\mathcal{G}}(c,U)_{(0)}=U$ is an
$A(V_{\mathcal{G}}(c,0))$-module on which $[v_X]$
acts as $X(0)$, hence by left multiplication by $X_0$
for $X=L,I,J,H$.
Thus, for every $u\in U(\fg_{\mathcal{G}})$,
$\bar\psi(u)$ acts on $U$ by left multiplication by $u$.
If $u\in\ker\bar\psi$, then
\begin{equation*}
	0=\bar\psi(u)\cdot1=u.
\end{equation*}
Thus, $\bar\psi$ is injective.		
\end{proof}

To determine $A(V^0_{\mathcal{G}}(c,0))$,
we now compute the classes of the spanning vectors in \eqref{eq:I0-central-span}.

\begin{lemma}\label{lem:zhu-p-center}
	For $n\geq2$, $m\geq1$, and $u\in V_{\mathcal{G}}(c,0)$, the following
	relations hold in $A(V_{\mathcal{G}}(c,0))$:
	\begin{align}
		[(L(-n)^p-\delta_{p\mid n}L(-pn))u]
		&=(-1)^{np}(n-1)[u]
		\big([v_L]^p-[v_L]\big),\label{eq:zhu-L-p}\\
		[X(-n)^pu]
		&=(-1)^n(n-1)[u][v_X]^p,\label{eq:zhu-IJ-p}\\
		[(H(-m)^p-H(-pm))u]
		&=(-1)^{m-1}[u]
		\big([v_H]^p-[v_H]\big),\label{eq:zhu-H-p}
	\end{align}
	where $X=I,J$.
\end{lemma}

\begin{proof}
	By linearity, we may assume that $u$ is homogeneous of degree $d$.
	We first consider \eqref{eq:zhu-L-p}.
	By \eqref{eq:u*v} and \eqref{eq:zhu-comm},
	\begin{equation*}
		[(L(-2)+L(-1))u]=[u][v_L].
	\end{equation*}
Using the same reduction as in
		\cite[Lemma~4.1]{JLM} with Virasoro generators $-L_k$,
		together with $L(0)u=-du$ and
		$[(L(-2)+L(-1))u]=[u][v_L]$, we obtain
	\begin{align}
		[L(-n)u]
		&=(-1)^n
		[\big((n-1)(L(-2)+L(-1))+L(0)\big)u]
		\nonumber\\
		&=(-1)^n[u]\big((n-1)[v_L]-d[\1]\big)
		\label{eq:zhu-L-n}
	\end{align}
	for $n\ge2$.
	Since $L(-n)^ju$ is homogeneous of degree $d+jn$,
	repeated application of \eqref{eq:zhu-L-n} gives
	\begin{equation}\label{eq:L-n-p-product}
		[L(-n)^pu]
		=(-1)^{np}[u]
		\prod_{j=0}^{p-1}
		\big((n-1)[v_L]-(d+jn)[\1]\big).
	\end{equation}
By \eqref{eq:app-fermat-x}, we have
	$\prod_{j=0}^{p-1}(x-jn)=x^p-n^{p-1}x$
	in $\BF[x]$.
	It follows from \eqref{eq:L-n-p-product}
	and Fermat's little theorem that
	\begin{align}
		[L(-n)^pu]
		&=(-1)^{np}[u]
		\big(
		\big((n-1)[v_L]-d[\1]\big)^p
		-n^{p-1}\big((n-1)[v_L]-d[\1]\big)
		\big)
		\nonumber\\
		&=(-1)^{np}[u]
		\big(
		(n-1)[v_L]^p-d[\1]
		-(1-\delta_{p\mid n})
		\big((n-1)[v_L]-d[\1]\big)
		\big)
		\nonumber\\
		&=(-1)^{np}(n-1)[u]
		\big([v_L]^p-[v_L]\big)
		-\delta_{p\mid n}(-1)^{np}[u]
		\big([v_L]+d[\1]\big).
		\label{eq:L-P-1}
	\end{align}
	On the other hand, applying \eqref{eq:zhu-L-n}
	with $n$ replaced by $pn$ gives
	\begin{equation}\label{eq:L-P-2}
		[L(-pn)u]
		=-(-1)^{np}[u]\big([v_L]+d[\1]\big).
	\end{equation}
	Thus, \eqref{eq:zhu-L-p} follows from
	\eqref{eq:L-P-1} and \eqref{eq:L-P-2}.
	
	We now consider \eqref{eq:zhu-IJ-p} for $X=I,J$.
	We claim that
	\begin{equation}\label{eq:zhu-X-n}
		[X(-n)u]
		=(-1)^n\big((n-1)[u][v_X]-[X(-1)u]\big)
	\end{equation}
	for $n\ge1$.
	The case $n=1$ is immediate.
	Since $\deg v_X=2$,
	\eqref{eq:u*v} and \eqref{eq:zhu-comm} give
	\begin{align*}
		[v_X][u]
		&=[X(-2)u]+2[X(-1)u]+[X(0)u],\\
		[v_X][u]-[u][v_X]
		&=[X(-1)u]+[X(0)u],
	\end{align*}
so \eqref{eq:zhu-X-n} holds for $n=2$.
	For $n\ge3$,
	\eqref{eq:u-circ-v} and \eqref{eq:def-O(V)} give
	\begin{equation*}
		v_X\circ_{n-3}u
		=X(-n)u+2X(-n+1)u+X(-n+2)u
		\in O(V_{\mathcal{G}}(c,0)).
	\end{equation*}
		Thus, \eqref{eq:zhu-X-n} holds for all $n\ge1$ by induction.

	Since $X(-n)$ and $X(-1)$ commute and
	$\binom pj=0$ in $\BF$ for $0<j<p$,
	iterating \eqref{eq:zhu-X-n} gives
	\begin{equation*}
		[X(-n)^pu]
		=(-1)^{np}
		\big((n-1)^p[u][v_X]^p-[X(-1)^pu]\big).
	\end{equation*}
	Since $X_{-1}^p$ is central in $U(\mathcal{G})$
	(see Corollary \ref{coro:p-cent}) and $X(-1)\1=0$,
	$X(-1)^p$ acts trivially on $V_{\mathcal{G}}(c,0)$.
	It follows from Fermat's little theorem that
	\begin{equation*}
		[X(-n)^pu]=(-1)^n(n-1)[u][v_X]^p,
	\end{equation*}
	which proves \eqref{eq:zhu-IJ-p}.
	
Finally, 
since $\deg v_H=1$,
the argument in the proof of \cite[Lemma~6.2]{JLM} gives \eqref{eq:zhu-H-p}.
\end{proof}

We then obtain the following result.
\begin{theorem}\label{th:iso-u(g)-V0}
	The algebra isomorphism $\bar{\psi}$ from $U(\fg_{\mathcal{G}})$ to $A(V_{\mathcal{G}}(c,0))$ in Theorem \ref{th:A(V)} descends to an isomorphism from $\mathfrak{u}(\fg_{\mathcal{G}})$ to $A(V^{0}_{\mathcal{G}}(c,0))$.
\end{theorem}

\begin{proof}
Set
$[I^0]=\{[u]\mid u\in I^0\}
\subset A(V_{\mathcal{G}}(c,0))$.
By \cite[Theorem~3.3(5)]{DR},
$[I^0]$ is an ideal of $A(V_{\mathcal{G}}(c,0))$ and
\begin{equation*}
	A(V^0_{\mathcal{G}}(c,0))
	\cong A(V_{\mathcal{G}}(c,0))/[I^0].
\end{equation*}
By \eqref{eq:I0-central-span} and Lemma \ref{lem:zhu-p-center},
$[I^0]$ is exactly the ideal generated by
\begin{equation}\label{eq:gen-I-0}
	[v_L]^p-[v_L],\quad
	[v_I]^p,\quad
	[v_J]^p,\quad
	[v_H]^p-[v_H].
\end{equation}
Under $\bar\psi$, the elements in \eqref{eq:gen-I-0}
are respectively the images of
\begin{equation}\label{eq:gen-p-ideal}
	L_0^p-L_0,\quad I_0^p,\quad J_0^p,\quad H_0^p-H_0.
\end{equation}
By \eqref{eq:rest-stru-g},
the elements in \eqref{eq:gen-p-ideal}
generate the defining ideal of $\mathfrak u(\fg_{\mathcal{G}})$
in $U(\fg_{\mathcal{G}})$.
Hence $\bar\psi$ induces the required isomorphism.
\end{proof}

We next classify the irreducible
$U(\fg_{\mathcal{G}})$-modules.
By the universal property of $U(\fg_{\mathcal{G}})$,
we may identify $\fg_{\mathcal{G}}$-modules with unital
$U(\fg_{\mathcal{G}})$-modules.
\begin{lemma}\label{lem:zero-mode-finite}
	Let $U$ be an irreducible $\fg_{\mathcal{G}}$-module. There exist
	unique scalars $\alpha,\beta,\gamma,\delta\in\BF$ such that, for every
	$u\in U$,
\begin{equation}\label{eq:zero-mode-central-scalars}
	L_0u=\alpha u, \;\;	(H_0^p-H_0)u=\beta u,\;\;	I_0^pu=\gamma^p u,	\;\;J_0^pu=\delta^p u,\;\;	I_0J_0u=\gamma\delta u.
\end{equation}	
\end{lemma}

\begin{proof}
By \eqref{eq:non-bra-g} and Corollary \ref{coro:p-cent},
$L_0$, $H_0^p-H_0$, $I_0^p$, and $J_0^p$ are central in $U(\fg_{\mathcal{G}})$.
Furthermore, 
\eqref{eq:H-I;H-J} gives
	\begin{equation*}
		[H_0,I_0J_0]
		=[H_0,I_0]J_0+I_0[H_0,J_0]=0.
	\end{equation*}
Hence $I_0J_0$ is also central in  $U(\fg_{\mathcal{G}})$.
	
Since $\fg_{\mathcal{G}}$ is finite-dimensional,
$U$ is finite-dimensional by \cite[Chapter~5, Theorem~2.4]{SF}.
By Schur's lemma, all five elements $L_0$, $H_0^p-H_0$, $I_0^p$, $J_0^p$ and $I_0J_0$   act on $U$ as scalars.
Since every element of $\BF$ has a unique $p$th root,
the first four relations in
\eqref{eq:zero-mode-central-scalars} hold with unique
$\alpha,\beta,\gamma,\delta\in\BF$.
Finally, $(I_0J_0)^p=I_0^pJ_0^p$ gives the last relation
by the uniqueness of $p$th roots.
\end{proof}

Let $\BF_p$ denote the prime subfield of $\BF$.
We now give explicit realizations of all irreducible modules for $\fg_{\mathcal{G}}$.

\begin{theorem}\label{th:classification-g-G}
	The irreducible $\fg_{\mathcal{G}}$-modules are,
	up to isomorphism, precisely the modules in the following two families.
	\begin{enumerate}[(i)]
		\item For $\alpha,\lambda\in\BF$, let
		$\BF_{\alpha,\lambda}=\BF v$ be the one-dimensional
		$\fg_{\mathcal{G}}$-module with action given by
		\begin{equation}\label{eq:zero-mode-one-dimensional}
			L_0v=\alpha v,\quad H_0v=\lambda v,\quad I_0v=J_0v=0.
		\end{equation}
		The modules $\BF_{\alpha,\lambda}$ and $\BF_{\alpha',\lambda'}$
		are isomorphic if and only if $\alpha=\alpha'$ and $\lambda=\lambda'$.
		
		\item For $\alpha,\lambda,\gamma,\delta\in\BF$ with
		$(\gamma,\delta)\ne(0,0)$, let
		$E(\alpha,\lambda;\gamma,\delta)$ be the $p$-dimensional
		$\fg_{\mathcal{G}}$-module with basis
		$\{e_i\mid i\in\BF_p\}$ and action given by
	\begin{equation}\label{eq:zero-mode-cyclic-actions}
		L_0e_i=\alpha e_i,\quad H_0e_i=(\lambda+i)e_i,\quad
		I_0e_i=\gamma e_{i+1},\quad J_0e_i=\delta e_{i-1}\quad\text{for } i\in\BF_p.
	\end{equation}
		The modules $E(\alpha,\lambda;\gamma,\delta)$ and
		$E(\alpha',\lambda';\gamma',\delta')$ are isomorphic if and only if
		\begin{equation}\label{eq:zero-mode-isomorphism}
			\alpha=\alpha',\quad \gamma=\gamma',\quad \delta=\delta',
			\quad \lambda-\lambda'\in\BF_p.
		\end{equation}
	\end{enumerate}
\end{theorem}

\begin{proof}
A direct calculation using \eqref{eq:non-bra-g} shows that $\BF_{\alpha,\lambda}$ and $E(\alpha,\lambda;\gamma,\delta)$ are $\fg_{\mathcal{G}}$-modules.
We now prove irreducibility.
The one-dimensional modules $\BF_{\alpha,\lambda}$
are clearly irreducible.
	
For $E(\alpha,\lambda;\gamma,\delta)$,
let $W$ be a nonzero submodule.
Choose a nonzero vector $w=\sum\limits_{j\in\BF_p}c_je_j\in W$ and take $i\in\BF_p$ with $c_i\ne0$.
Since $H_0W\subset W$,
\eqref{eq:zero-mode-cyclic-actions} gives
\begin{equation*}
	\bigg(\prod_{\substack{j\in\BF_p\\j\ne i}}\frac{H_0-(\lambda+j)\operatorname{id}_W}{i-j}\bigg)w=c_ie_i\in W.
\end{equation*}
Hence $e_i\in W$.
If $\gamma\ne0$, 
then $I_0^ke_i=\gamma^ke_{i+k}\in W$ for $0\le k<p$.
If $\gamma=0$, 
then $\delta\ne0$ and $J_0^ke_i=\delta^ke_{i-k}\in W$ for $0\le k<p$.
In either case, $e_j\in W$ for every $j\in\BF_p$,
so $W=E(\alpha,\lambda;\gamma,\delta)$.
	
We next prove that every irreducible $\fg_{\mathcal{G}}$-module is isomorphic to $\BF_{\alpha,\lambda}$ or $E(\alpha,\lambda;\gamma,\delta)$.
Let $U$ be an irreducible $\fg_{\mathcal{G}}$-module with scalars $\alpha,\beta,\gamma,\delta$ as in Lemma \ref{lem:zero-mode-finite}.
Set
\begin{equation*}
	K_I=\{u\in U\mid I_0u=0\},\quad
	K_J=\{u\in U\mid J_0u=0\}.
\end{equation*}
By \eqref{eq:non-bra-g},
both $K_I$ and $K_J$ are $\fg_{\mathcal{G}}$-submodules.
By \eqref{eq:zero-mode-central-scalars},	
$I_0$ acts nilpotently if $\gamma=0$ and invertibly otherwise; 
the same holds for $J_0$ and $\delta$.
The irreducibility of $U$ therefore gives
\begin{equation*}
	K_I=
	\begin{cases}
		U,&\gamma=0,\\
		\{0\},&\gamma\ne0,
	\end{cases}
	\qquad
	K_J=
	\begin{cases}
		U,&\delta=0,\\
		\{0\},&\delta\ne0.
	\end{cases}
\end{equation*}
	
Since $U$ is finite-dimensional and $\BF$ is algebraically closed, 
we can choose a nonzero vector $e_0\in U$ with $H_0e_0=\lambda e_0$.
If $\gamma=\delta=0$, 
then $\BF e_0$ is a nonzero $\fg_{\mathcal{G}}$-submodule,
so $U=\BF e_0\cong\BF_{\alpha,\lambda}$.
	
Now suppose that $\gamma\ne0$.
Set
\begin{equation*}
	e_i:=\gamma^{-i}I_0^ie_0
	\quad\text{for }0\le i<p.
\end{equation*}
Using $[H_0,I_0]=I_0$ and $I_0^pe_0=\gamma^pe_0$ (see \eqref{eq:zero-mode-central-scalars}), 
we have
\begin{equation*}
	H_0e_i=(\lambda+i)e_i,\quad
	I_0e_i=\gamma e_{i+1}
	\quad\text{for }i\in\BF_p.
\end{equation*}
	Since $K_I=\{0\}$, every $e_i$ is nonzero.
The vectors $e_i$ for $i\in\BF_p$ have distinct
$H_0$-eigenvalues $\lambda+i$ and are thus linearly independent.
Moreover, 
\eqref{eq:zero-mode-central-scalars} gives
\begin{equation*}
	I_0(J_0e_i-\delta e_{i-1})
	=\gamma\delta e_i-\gamma\delta e_i=0
	\quad\text{for }i\in\BF_p.
\end{equation*}
Since $K_I=\{0\}$, 
we obtain $J_0e_i=\delta e_{i-1}$ for $i\in\BF_p$.
Thus,
$\bigoplus\limits_{i\in\BF_p}\BF e_i$ is a nonzero $\fg_{\mathcal{G}}$-submodule of $U$.
By irreducibility of $U$,
\begin{equation*}
	U=\bigoplus_{i\in\BF_p}\BF e_i
	\cong E(\alpha,\lambda;\gamma,\delta).
\end{equation*}
	
Finally, 
suppose that $\gamma=0$ and $\delta\ne0$,
so $I_0$ acts trivially on $U$ and $J_0$ acts invertibly on $U$.
Set
\begin{equation*}
	e_{-i}:=\delta^{-i}J_0^ie_0
	\quad\text{for }0\le i<p.
\end{equation*}
Using $[H_0,J_0]=-J_0$ and $J_0^pe_0=\delta^pe_0$ (see \eqref{eq:zero-mode-central-scalars}),
we obtain
\begin{equation*}
	H_0e_{-i}=(\lambda-i)e_{-i},\quad
	J_0e_{-i}=\delta e_{-i-1}
	\quad\text{for }0\le i<p.
\end{equation*}
As in the case $\gamma\ne0$,
the vectors $e_{-i}$ for $i\in\BF_p$ form a basis of $U$,
and hence $U\cong E(\alpha,\lambda;0,\delta)$.
	
We now determine when two listed modules are isomorphic.
By \eqref{eq:zero-mode-one-dimensional},
$\BF_{\alpha,\lambda}$ and $\BF_{\alpha',\lambda'}$ are isomorphic if and only if $\alpha=\alpha'$ and $\lambda=\lambda'$.
	
If $E(\alpha,\lambda;\gamma,\delta)$ and
$E(\alpha',\lambda';\gamma',\delta')$ are isomorphic,
comparing the actions of
$L_0$, $H_0^p-H_0$, $I_0^p$, and $J_0^p$ gives
\begin{equation*}
	\alpha=\alpha',\quad
	\lambda^p-\lambda=(\lambda')^p-\lambda',\quad
	\gamma^p=(\gamma')^p,\quad
	\delta^p=(\delta')^p.
\end{equation*}
The injectivity of the Frobenius map gives
$\gamma=\gamma'$ and $\delta=\delta'$.
Moreover, $(\lambda-\lambda')^p=\lambda-\lambda'$,
so $\lambda-\lambda'\in\BF_p$.
Hence \eqref{eq:zero-mode-isomorphism} holds.
	
Conversely, suppose that \eqref{eq:zero-mode-isomorphism}
holds and write $\lambda'=\lambda+r$ with $r\in\BF_p$.
The assignment $f(e_i)=e'_{i-r}$ for $i\in\BF_p$
defines a linear isomorphism $f$ from
$E(\alpha,\lambda;\gamma,\delta)$ to
$E(\alpha',\lambda';\gamma',\delta')$.
The actions in \eqref{eq:zero-mode-cyclic-actions} give
\begin{equation*}
	f(X_0e_i)=X_0f(e_i)
	\quad\text{for }X_0\in\{L_0,H_0,I_0,J_0\}
	\text{ and }i\in\BF_p.
\end{equation*}
Hence $f$ is a $\fg_{\mathcal{G}}$-module isomorphism.
\end{proof}

We next classify the irreducible modules for the restricted enveloping algebra $\mathfrak u(\fg_{\mathcal{G}})$.

\begin{corollary}\label{coro:zero-mode-restricted}
The $p^2$ one-dimensional modules $\BF_{\alpha,\lambda}$ with $(\alpha,\lambda)\in\BF_p^2$
form a complete set of pairwise nonisomorphic irreducible $\mathfrak{u}(\fg_{\mathcal{G}})$-modules.
\end{corollary}

\begin{proof}
By the definition of $\mathfrak{u}(\fg_{\mathcal{G}})$ and \eqref{eq:rest-stru-g},
the action of $\fg_{\mathcal{G}}$ on a module $U$ factors through $\mathfrak{u}(\fg_{\mathcal{G}})$ if and only if
\begin{equation*}
	L_0^pu=L_0u,\quad H_0^pu=H_0u,\quad I_0^pu=J_0^pu=0
	\quad\text{for }u\in U.
\end{equation*}
By \eqref{eq:zero-mode-one-dimensional},
the action on $\BF_{\alpha,\lambda}$ factors through $\mathfrak u(\fg_{\mathcal{G}})$ if and only if $\alpha^p=\alpha$ and $\lambda^p=\lambda$,
equivalently, $\alpha,\lambda\in\BF_p$.
	
By \eqref{eq:zero-mode-cyclic-actions},
$I_0^p$ and $J_0^p$ act on
$E(\alpha,\lambda;\gamma,\delta)$
as the scalars $\gamma^p$ and $\delta^p$, respectively.	
Since $(\gamma,\delta)\ne(0,0)$,
at least one of $\gamma^p$ and $\delta^p$ is nonzero, 
so the action does not factor through $\mathfrak{u}(\fg_{\mathcal{G}})$.
The result now follows from Theorem \ref{th:classification-g-G}.
\end{proof}

Combining Theorems \ref{th:A(V)}, \ref{th:classification-g-G},
and \cite[Propositions~2.10, 2.11]{JLM},
we  obtain the following classification.
\begin{theorem}
The following modules form a complete set of pairwise nonequivalent irreducible $\BN$-graded $V_{\mathcal{G}}(c,0)$-modules:
	\begin{enumerate}[(i)]
	\item $L_{\mathcal{G}}(c,\BF_{\alpha,\lambda})$
	for $\alpha,\lambda\in\BF$.
	
	\item $L_{\mathcal{G}}
(c,E(\alpha,\lambda;\gamma,\delta))$
for $\alpha,\gamma,\delta\in\BF$
with $(\gamma,\delta)\ne(0,0)$,
where $\lambda$ ranges over a fixed set of representatives
for the additive cosets of $\BF_p$ in $\BF$.
\end{enumerate}
\end{theorem}
	
Furthermore, 
Theorem \ref{th:iso-u(g)-V0}, Corollary \ref{coro:zero-mode-restricted}, and \cite[Propositions 2.10, 2.11]{JLM}
give the following classification.

\begin{theorem}\label{th:irr-V0-mod}
	The vertex algebra $V^0_{\mathcal{G}}(c,0)$ has exactly $p^2$
	irreducible $\BN$-graded modules up to equivalence, namely
	\begin{equation*}
		L_{\mathcal{G}}(c,\BF_{\alpha,\lambda})
		\quad\text{for }\alpha,\lambda\in\BF_p.
	\end{equation*}
\end{theorem}

\section{Invariant bilinear forms and the simple graded quotient of $V_{\mathcal{G}}(c,0)$}\label{sec:5}
In this section, we determine the space of invariant bilinear
forms on $V_{\mathcal{G}}(c,0)$,
where $c=(c_1,c_2,c_3)\in\BF^3$.
For $c_2=0$, we identify the simple graded quotient
and classify the irreducible $\BN$-graded modules
of the quotient for certain values of $c$.

We first recall the (cocommutative) bialgebra $\mathcal{H}$
(see \cite{LM2}).
As an associative algebra over $\BF$,
$\mathcal{H}$ is generated by
\begin{equation*}
	L_{(-1)}^{(k)},\quad L_{(0)}^{(k)},\quad L_{(1)}^{(k)}
	\qquad\text{for }k\in\BN,
\end{equation*}
subject to $L_{(-1)}^{(0)}=L_{(0)}^{(0)}=L_{(1)}^{(0)}=1$
and the following relations for $j,k\in\BN$:
\begin{align}
	L_{(0)}^{(j)}L_{(0)}^{(k)}
	&=\sum_{t=0}^{j}\binom{j}{t}\binom{k+t}{j}
	L_{(0)}^{(k+t)},\label{eq:alg-str}\\
	L_{(\pm1)}^{(j)}L_{(\pm1)}^{(k)}
	&=\binom{j+k}{j}L_{(\pm1)}^{(j+k)},\\
	L_{(0)}^{(j)}L_{(\pm1)}^{(k)}
	&=\sum_{t=0}^{j}\binom{\pm2k}{t}
	L_{(\pm1)}^{(k)}L_{(0)}^{(j-t)},\\
	L_{(1)}^{(j)}L_{(-1)}^{(k)}
	&=\sum_{s=0}^{\min(j,k)}\sum_{t=0}^{s}
	\binom{-j-k+2s}{t}(-1)^s
	L_{(-1)}^{(k-s)}L_{(0)}^{(s-t)}L_{(1)}^{(j-s)}.
\end{align}
Its coalgebra structure is given by
\begin{equation}\label{eq:bialg-str}
	\varepsilon(L_{(r)}^{(k)})=\delta_{k,0},
	\quad
	\Delta(L_{(r)}^{(k)})
	=\sum_{i=0}^{k}
	L_{(r)}^{(k-i)}\otimes L_{(r)}^{(i)}
\end{equation}
for $r\in\{-1,0,1\}$ and $k\in\BN$.

Note that $\mathcal{H}$ is a $\BZ$-graded algebra with
\begin{equation}
	\deg L_{(\pm1)}^{(k)}=\mp k,\quad
	\deg L_{(0)}^{(k)}=0
	\quad\text{for }k\in\BN.
\end{equation}
A $\BZ$-graded $\mathcal{H}$-module is defined in the usual way.
Furthermore, recall from \cite[Definition 3.2]{LM2} that
a {\em $\BZ$-graded weight $\mathcal{H}$-module} is a $\BZ$-graded
$\mathcal{H}$-module $W=\bigoplus\limits_{m\in\BZ} W_{(m)}$
on which $L_{(0)}^{(k)}$
acts as $\binom{-2 \deg }{k}$ for $k\in \BN$, namely
\begin{equation}
	L_{(0)}^{(k)}|_{W_{(m)}}=\binom{-2m}{k}
	\;\text{ for }k\in \BN, m\in \BZ.
\end{equation}
	We use the formal generating function
\begin{equation*}
		e^{xL_{(\pm 1)}}:=\sum_{k\in\BN}x^kL_{(\pm 1)}^{(k)}.
\end{equation*}

We then recall the notion of an $\mathcal{H}$-module vertex algebra (over $\BF$),
 which  is an analogue of a quasi-vertex operator algebra in characteristic zero (see \cite[\S 2.8]{FHL}).
\begin{definition}[{\cite[Definition~3.7]{LM2}}]
	\label{def:H-mod-va}
	An {\em $\mathcal{H}$-module vertex algebra}
	is a $\BZ$-graded vertex algebra
	$V=\bigoplus_{m\in\BZ}V_{(m)}$
	which is also a $\BZ$-graded weight $\mathcal{H}$-module,
		with $L_{(-1)}^{(k)}=\D^{(k)}$ on $V$ for $k\in\BN$,
	satisfying the following conditions:
	\begin{enumerate}[(i)]
		\item $V_{(m)}=0$ for $m$ sufficiently negative.
		\item 
			$L_{(1)}^{(k)}\1=\delta_{k,0}\1$
		for $k\in\BN$.
		\item For homogeneous $v\in V$,
		\begin{equation}\label{eq:(iii)-def-H-mod-VA}
			e^{xL_{(1)}}Y(v,x_0)e^{-xL_{(1)}}
			=Y\big(
			e^{x(1-xx_0)L_{(1)}}
			(1-xx_0)^{-2\deg}v,
			x_0/(1-xx_0)\big).
		\end{equation}
	\end{enumerate}
\end{definition}

We next show that $V_{\mathcal{G}}(c,0)$ carries an
$\mathcal{H}$-module vertex algebra structure.
For $i=1,2,3$, $X=L,I,J$, $k\in\BN$, and $m\in\BZ$, define
\begin{align}
	L_{(-1)}^{(k)}\bc_i
	&=\delta_{k,0}\bc_i,
	&
	L_{(0)}^{(k)}\bc_i
	&=\delta_{k,0}\bc_i,
	&
	L_{(1)}^{(k)}\bc_i
	&=\delta_{k,0}\bc_i;
	\label{eq:H-mod-acts-1}\\
	L_{(-1)}^{(k)}X_m
	&=(-1)^k\binom{m+1}{k}X_{m-k},
	&
	L_{(0)}^{(k)}X_m
	&=\binom{2m}{k}X_m,
	&
	L_{(1)}^{(k)}X_m
	&=\binom{-m+1}{k}X_{m+k};
	\label{eq:H-mod-acts-2}\\
	L_{(-1)}^{(k)}H_m
	&=(-1)^k\binom{m}{k}H_{m-k},
	&
	L_{(0)}^{(k)}H_m
	&=\binom{2m}{k}H_m,
	&
	L_{(1)}^{(k)}H_m
	&=\binom{-m}{k}H_{m+k}-2(1-\delta_{k,0})\nonumber\\
	&&&&&\?\times\delta_{m+k,0}\bc_2.
	\label{eq:H-mod-acts-3}
\end{align}

\begin{lemma}\label{lem:H-module-Galilean}
	The Lie algebra $\mathcal{G}$ is an $\mathcal{H}$-module Lie algebra with the
action given by \eqref{eq:H-mod-acts-1}--\eqref{eq:H-mod-acts-3}.
\end{lemma}

\begin{proof}
For $X=L,I,J$, the actions in
\eqref{eq:H-mod-acts-2} make
$\bigoplus_{m\in\BZ}\BF X_m$ an $\mathcal{H}$-module
by \cite[(5.12)]{LM2}.
By \cite[Example~3.4]{LM2}, both $\BF[x]$ and
$\BF[x^{-1}]$ are $\mathcal{H}$-submodules of $\BF[x,x^{-1}]$.
The actions in \eqref{eq:H-mod-acts-3} give
$\mathcal{H}$-module isomorphisms
\begin{equation*}
	\BF[x]\cong\bigoplus_{m\ge0}\BF H_m,\quad
	\BF[x^{-1}]\cong
	\big(\bigoplus_{m<0}\BF H_m\big)\oplus\BF(H_0-2\bc_2),
\end{equation*}
under $x^m\mapsto H_m$ for $m\ge0$ and $x^m\mapsto H_m-2\delta_{m,0}\bc_2$ for $m\le0$,
respectively.
Since $p>2$, we have
\begin{equation*}
	\Big(\bigoplus_{m\ge0}\BF H_m\Big)
	\oplus
	\Big(\bigoplus_{m<0}\BF H_m\Big)
	\oplus\BF(H_0-2\bc_2)
	=
	\Big(\bigoplus_{m\in\BZ}\BF H_m\Big)\oplus\BF\bc_2.
\end{equation*}
Hence
$\big(\bigoplus_{m\in\BZ}\BF H_m\big)\oplus\BF\bc_2$
is an $\mathcal{H}$-module.
Together with the trivial $\mathcal{H}$-modules
$\BF\bc_1$ and $\BF\bc_3$, this proves that
$\mathcal{G}$ is an $\mathcal{H}$-module.

Since $L_{(-1)}^{(k)}=\D^{(k)}$ on $\mathcal{G}$
for $k\in \BN$ (see \eqref{eq:B-mod-acts-G}
and \eqref{eq:H-mod-acts-1}--\eqref{eq:H-mod-acts-3}),
it follows from
Lemma \ref{lem:U(G)-B-mod-alg} that
\begin{equation*}
	L_{(-1)}^{(k)}[a,b]=\D^{(k)}[a,b]=\sum_{i=0}^{k}
	[\D^{(k-i)}a,\D^{(i)}b]
	=\sum_{i=0}^{k}
	[L_{(-1)}^{(k-i)}a,L_{(-1)}^{(i)}b]
\end{equation*}
for $k\in \BN$ and $a,b\in \mathcal{G}$.

We next verify that
\begin{equation}\label{eq:H-coproduct-Galilean}
	L_{(1)}^{(k)}[a,b]
	=\sum_{i=0}^{k}
	[L_{(1)}^{(k-i)}a,L_{(1)}^{(i)}b]
\end{equation}
for $k\in\BN$ and $a,b\in\mathcal{G}$.
The case $k=0$ is immediate.
For $k\in\BZ_+$, by bilinearity and skew-symmetry,
it suffices to verify \eqref{eq:H-coproduct-Galilean}
for the pairs of basis elements appearing in \eqref{eq:Lm-Ln-brac}--\eqref{eq:H-I;H-J},
as the remaining cases are immediate.

For $a=L_m$ and $b=L_n$,
\eqref{eq:H-coproduct-Galilean} follows from
\cite[(5.11)--(5.13)]{LM2},
with $-L_n$ taken as the Virasoro generators.
The calculation for the noncentral part also proves \eqref{eq:H-coproduct-Galilean}
for $a=L_m$ and $b=X_n$ with $X=I,J$.

For $a=H_m$ and $b=H_n$, the calculation in
\cite[Lemma~5.1]{LM2} applies, since the additional terms in
$L_{(1)}^{(k)}H_m$ are central.
For $a=H_m$ and $b=X_n$ with $X=I,J$,
	\eqref{eq:H-coproduct-Galilean} follows from
	\eqref{eq:H-I;H-J} and \eqref{eq:Vand-iden}.

It remains to consider $a=L_m$ and $b=H_n$.
By \eqref{eq:H-mod-acts-1}--\eqref{eq:H-mod-acts-3},
\eqref{eq:Hm-Hn-brac}, and \eqref{eq:Vand-iden},
 we obtain
\begin{align*}
	\sum_{i=0}^{k}
	[L_{(1)}^{(k-i)}L_m,L_{(1)}^{(i)}H_n]
	&=
	\sum_{i=0}^{k}
	\binom{-m+1}{k-i}\binom{-n}{i}
	[L_{m+k-i},H_{n+i}]\\
	&=
	\sum_{i=0}^{k}
	\binom{-m+1}{k-i}\binom{-n}{i}(n+i)
	H_{m+n+k}\\
	&\?+
	\delta_{m+n+k,0}
	\sum_{i=0}^{k}
	\binom{-m+1}{k-i}\binom{-n}{i}
	(m+k-i)(m+k-i+1)\bc_2\\
	&=	n\sum_{i=0}^{k}
	\binom{-m+1}{k-i}\binom{-n-1}{i}
	H_{m+n+k}\\
	&\?+
	\delta_{m+n+k,0}
	\sum_{i=0}^{k}
	\binom{-m+1}{k-i}\binom{-n}{i}
	(n+i)(n+i-1)\bc_2\\
	&=	n\binom{-m-n}{k}H_{m+n+k}\\
	&\?+\delta_{m+n+k,0}
	\sum_{i=0}^{k}
	\binom{-m+1}{k-i}\bigg(	n(n+1)\binom{-n-2}{i}
	-2n\binom{-n-1}{i}\bigg)\bc_2\\
	&=
	n\binom{-m-n}{k}H_{m+n+k}
	-2n\delta_{m+n+k,0}\bc_2\\
	&=L_{(1)}^{(k)}[L_m,H_n].
\end{align*}
Here the penultimate equality uses $k\in \BZ_+$.
This proves \eqref{eq:H-coproduct-Galilean}.
Since $\mathcal{H}$ is generated by $L_{(\pm1)}^{(k)}$ for $k\in\BN$ (see \cite[Remark~4.2]{LM2}),
$\mathcal{G}$ is an $\mathcal{H}$-module Lie algebra.
\end{proof}

By Lemmas \ref{lem:H-module-Galilean}
and \ref{lem:H-mod-Lie-H-mod-alg},
$U(\mathcal{G})$ is an $\mathcal{H}$-module algebra.
The actions in \eqref{eq:H-mod-acts-1}--\eqref{eq:H-mod-acts-3}
preserve $\mathfrak p+\sum_{i=1}^3\BF(\bc_i-c_i)$
and hence the left ideal $J_c$ generated by this subspace.
Thus $V_{\mathcal{G}}(c,0)\cong U(\mathcal{G})/J_c$
inherits a $(U(\mathcal{G}),\mathcal{H})$-module structure.

Moreover, the actions in
\eqref{eq:H-mod-acts-1}--\eqref{eq:H-mod-acts-3}
are homogeneous with respect to the grading
\eqref{eq:deg-bas-G},
so $U(\mathcal{G})$ is a $\BZ$-graded $\mathcal{H}$-module.
Since $J_c$ is generated by homogeneous elements,
$J_c$ is a graded $\mathcal{H}$-submodule.
Hence $V_{\mathcal{G}}(c,0)$, equipped with the grading
\eqref{eq:Galilean-VOA-grading}, inherits a
$\BZ$-graded $\mathcal{H}$-module structure.
We then obtain the following result.
\begin{proposition}\label{prop:H-module-Galilean-VA}
	For $c=(c_1,c_2,c_3)\in\BF^3$, the vertex algebra
	$V_{\mathcal{G}}(c,0)$ is an $\mathcal{H}$-module vertex algebra.
\end{proposition}

\begin{proof}
As in \cite[Proposition~5.2]{LM2}, we have
$L_{(-1)}^{(k)}=\D^{(k)}$ and
$L_{(0)}^{(k)}=\binom{-2\deg}{k}$
on $V_{\mathcal{G}}(c,0)$ for $k\in\BN$.
Thus, $V_{\mathcal{G}}(c,0)$ is a $\BZ$-graded weight $\mathcal{H}$-module with the natural $\B$-module action.
Furthermore, $L_{(1)}^{(k)}\1=\delta_{k,0}\1$ for $k\in\BN$.
Since $V_{\mathcal{G}}(c,0)$ is $\BN$-graded,
conditions (i) and (ii) of Definition \ref{def:H-mod-va} hold.
It remains to verify condition (iii).

	By the same argument as in \cite[Proposition~5.5]{LM2},
Definition \ref{def:H-mod-va}(iii) holds for $v=v_L,v_I,v_J$.
	For $v=v_H$, we have $\deg v_H=1$.
	By \eqref{eq:H-mod-acts-3}
	and $L_{(1)}^{(k)}\1=\delta_{k,0}\1$,
	\begin{equation}\label{eq:L-1-vH}
		L_{(1)}^{(1)}v_H=-2c_2\1,\qquad
		L_{(1)}^{(k)}v_H=0
		\quad\text{for }k\ge2.
	\end{equation}
		Using \eqref{eq:H-mod-acts-3}
		and the $(U(\mathcal{G}),\mathcal{H})$-module structure, we obtain
	\begin{align*}
		&e^{xL_{(1)}}Y(v_H,x_0)e^{-xL_{(1)}}\\
		&=
		\sum_{n\in\BZ}\sum_{k\in\BN}
		\bigg(
		\binom{-n}{k}H(n+k)
		-2(1-\delta_{k,0})\delta_{n+k,0}c_2
		\bigg)x^k x_0^{-n-1}\\
		&=
		\sum_{m\in\BZ}H(m)x_0^{-m-1}
		\sum_{k\in\BN}\binom{k-m}{k}(xx_0)^k
		-2c_2\sum_{k\in\BZ_+}x^k x_0^{k-1}\\
		&=
		(1-xx_0)^{-2}
		Y\big(v_H,x_0/(1-xx_0)\big)
		-\frac{2c_2x}{1-xx_0}
		\operatorname{id}_{V_{\mathcal{G}}(c,0)}\\
		&=
		Y\big(
		e^{x(1-xx_0)L_{(1)}}(1-xx_0)^{-2\deg}v_H,
		x_0/(1-xx_0)\big).
	\end{align*}
	The last equality follows from \eqref{eq:L-1-vH}
	and $Y(\1,x_0)=\operatorname{id}_{V_{\mathcal{G}}(c,0)}$.
	Thus, Definition \ref{def:H-mod-va}(iii) holds
	for the generators $v_L,v_I,v_J,v_H$,
	and hence for all $v\in V_{\mathcal{G}}(c,0)$
	by \cite[Lemma~4.13]{LM2}.
\end{proof}

If $c_2=0$, then \eqref{eq:H-mod-acts-3} gives
\begin{equation}\label{eq:L-1-k-Hm}
	L_{(1)}^{(k)}(H_m)=\binom{-m}{k}H_{m+k}
\end{equation}
on $V_{\mathcal{G}}(c,0)$ for $m\in\BZ$ and $k\in\BN$.
Using \eqref{eq:L-1-k-Hm},
the coproduct formula \eqref{eq:bialg-str},
and the PBW basis \eqref{eq:Galilean-VOA-PBW},
we obtain the following lemma by adapting the arguments in
\cite[Proposition~5.3 and Lemma~5.6]{LM2}.

\begin{lemma}\label{lem:Galilean-H-positive-zero}
Let $c=(c_1,c_2,c_3)\in\BF^3$.
If $c_2=0$, then
$L_{(1)}^{(n)}V_{\mathcal{G}}(c,0)_{(n)}=0$
for all $n\in \BZ_+$.
\end{lemma}

Recall from \cite[Definition 4.1, Remark 4.2]{LM2}  that a bilinear form
$(\cdot,\cdot)$ on $V_{\mathcal{G}}(c,0)$ is {\em invariant} if
\begin{align}
\big(Y(u,x)v,w\big)&=\big(v,Y(e^{xL_{(1)}}(-x^{-2})^{\deg}u,x^{-1})w\big),\label{eq:Galilean-form-VA-inv}\\
(e^{xL_{(\pm 1)}}v,w)&=(v,e^{xL_{(\mp 1)}}w)\label{eq:Galilean-form-H}
\end{align}
for $u,v,w\in V_{\mathcal{G}}(c,0)$.

Since $V_{\mathcal{G}}(c,0)$ is $\BN$-graded with $V_{\mathcal{G}}(c,0)_{(0)}=\BF\1$,
Proposition \ref{prop:H-module-Galilean-VA},	Lemma \ref{lem:Galilean-H-positive-zero}, 
and \cite[Corollary~4.11]{LM2} give the following result.
\begin{theorem}\label{th:Galilean-invariant-forms}
Let $c=(c_1,c_2,c_3)\in\BF^3$ with $c_2=0$.
Then the space of invariant bilinear forms on $V_{\mathcal{G}}(c,0)$ is one-dimensional 
and there exists a unique symmetric invariant bilinear form $(\cdot,\cdot)$ on $V_{\mathcal{G}}(c,0)$ such that $(\1,\1)=1$.
\end{theorem}

\begin{remark}\label{rmk:c2 ne 0}
If $c_2\ne0$, then the space of invariant bilinear forms on $V_{\mathcal{G}}(c,0)$ is zero.
Indeed, since $p>2$,
\eqref{eq:L-1-vH} gives $L_{(1)}^{(1)}v_H=-2c_2\1\ne0$,
so
\begin{equation*}
	L_{(1)}^{(1)}V_{\mathcal{G}}(c,0)_{(1)}
	=\BF\1=V_{\mathcal{G}}(c,0)_{(0)}.
\end{equation*}
The conclusion follows from \cite[Theorem~4.10]{LM2}.
\end{remark}

For the remainder of this section, 
assume $c_2=0$ and let
$(\cdot,\cdot)$ denote the unique symmetric invariant bilinear
form on $V_{\mathcal{G}}(c,0)$ such that $(\1,\1)=1$
(see Theorem~\ref{th:Galilean-invariant-forms}).
By invariance (see  \eqref{eq:Galilean-form-VA-inv}), we have
\begin{equation}\label{eq:inv}
	(X(m)u,v)=(u,X(-m)v),\qquad
	(H(m)u,v)=-(u,H(-m)v)
\end{equation}
for $u,v\in V_{\mathcal{G}}(c,0)$, $m\in\BZ$, and $X=L,I,J$.
Set
\begin{equation*}
	\Rad(\cdot,\cdot)
	:=\{u\in V_{\mathcal{G}}(c,0)\mid
	(u,v)=0\text{ for all }v\in V_{\mathcal{G}}(c,0)\}.
\end{equation*}

\begin{lemma}\label{lem:rad-max-ideal}
$\Rad(\cdot,\cdot)$ is the unique maximal nontrivial
graded ideal of $V_{\mathcal{G}}(c,0)$.
\end{lemma}

\begin{proof}
Since $V_{\mathcal{G}}(c,0)$ is an $\BN$-graded
		$\mathcal{H}$-module vertex algebra by Proposition \ref{prop:H-module-Galilean-VA},
		its adjoint module is a lower truncated
		$(V_{\mathcal{G}}(c,0),\mathcal{H})$-module.
		By \cite[Lemma~4.3]{LM2},
	\begin{equation}\label{eq:orth-m-n}
		\big(V_{\mathcal{G}}(c,0)_{(m)},
		V_{\mathcal{G}}(c,0)_{(n)}\big)=0
		\quad\text{for }m,n\in\BN\text{ with }m\ne n.
	\end{equation}
Consequently,  $\Rad(\cdot,\cdot)$ is graded.
By \eqref{eq:Galilean-form-VA-inv} and
\eqref{eq:Galilean-form-H}, $\Rad(\cdot,\cdot)$ is stable
under $u_n$ and $\D^{(r)}=L_{(-1)}^{(r)}$
for $u\in V_{\mathcal{G}}(c,0)$, $n\in\BZ$, and $r\in\BN$.
		Hence $\Rad(\cdot,\cdot)$ is an ideal,
		which is proper since $(\1,\1)=1$.

Let $J$ be any proper graded ideal of $V_{\mathcal{G}}(c,0)$.
Since $\1\notin J$ and $V_{\mathcal{G}}(c,0)_{(0)}=\BF\1$, 
we have $J_{(0)}=0$.
Thus, $(\1,J)=0$ by \eqref{eq:orth-m-n}.
For homogeneous $u\in V_{\mathcal{G}}(c,0)$
and  $v\in J$,
\eqref{eq:Galilean-form-VA-inv} and $u_{-1}\1=u$ give
\begin{align*}
	(u,v)
	&=\Res_x x^{-1}\big(Y(u,x)\1,v\big)\\
	&=\Res_x x^{-1}
	\big(\1,Y(e^{xL_{(1)}}(-x^{-2})^{\deg u}u,x^{-1})v\big)
	=0,
\end{align*}
since $J$ is an ideal and $(\1,J)=0$.
By linearity and symmetry,
$J\subset \Rad(\cdot,\cdot)$.

In particular, 
the nonzero proper graded ideal $I^0$ is contained in $\Rad(\cdot,\cdot)$,
so $\Rad(\cdot,\cdot)\ne0$.
Therefore, $\Rad(\cdot,\cdot)$ is the unique maximal nontrivial graded ideal of $V_{\mathcal{G}}(c,0)$.
\end{proof}

It follows that
\begin{equation}
	L_{\mathcal{G}}(c,0)
	:=V_{\mathcal{G}}(c,0)/\Rad(\cdot,\cdot)
\end{equation}
is the unique simple graded quotient of $V_{\mathcal{G}}(c,0)$.

\begin{notation}
	\begin{enumerate}[(i)]
		\item Let $\mathcal L$ be the subalgebra of $\mathcal{G}$
		with basis
		$\{-L_n,H_n\mid n\in\BZ\}\cup\{\bc_1,\bc_2,\bc_3\}$.
		By \eqref{eq:Lm-Ln-brac} and \eqref{eq:Hm-Hn-brac},
		$\mathcal L$ is
			the modular analogue over $\BF$ of the rank-one twisted
			Heisenberg--Virasoro Lie algebra
		in \cite[Definition~2.4]{GW}.
		We denote by $V_{\mathcal L}(c,0)$
		the corresponding universal twisted
		Heisenberg--Virasoro vertex algebra.
		
		\item Let $\Vir$ be the Virasoro subalgebra of $\mathcal{G}$
		with basis
		$\{-L_n\mid n\in\BZ\}\cup\{\bc_1\}$.
		We denote by $V_{\Vir}(c_1,0)$ the universal
		Virasoro vertex algebra of central charge $c_1$,
		and by $L_{\Vir}(c_1,0)$ its simple graded quotient.
		
		\item Let $\fh=\BF H$ be an abelian Lie algebra
		with symmetric bilinear form $\langle H,H\rangle=1$.
			By \eqref{eq:Hm-Hn-brac}, we identify its affine Lie algebra
			$\widehat{\fh}$ with the subalgebra of $\mathcal{G}$
		with basis $\{H_n\mid n\in\BZ\}\cup\{\bc_3\}$.
		We denote by $V_{\widehat{\fh}}(c_3,0)$
		the corresponding Heisenberg vertex algebra of level $c_3$,
		and by $L_{\widehat{\fh}}(c_3,0)$
		its simple graded quotient.
	\end{enumerate}
\end{notation}

	For $c_3\ne0$, the proof of \cite[Theorem~3.16]{GW}
	applies over $\BF$, using the Heisenberg conformal vector
	$(2c_3)^{-1}H(-1)^2\1$ in \cite[Lemma~5.12]{LM}.
\begin{lemma}[{\cite{GW}}]\label{lem:HV-tensor}
	Let $c=(c_1,0,c_3)\in\BF^3$ with $c_3\ne0$.
	Then $V_{\mathcal L}(c,0)$ is isomorphic to
	$V_{\widehat{\fh}}(c_3,0)\otimes
	V_{\Vir}(c_1-1,0)$
	as an $\BN$-graded vertex algebra.
\end{lemma}

We now identify the simple graded quotient
$L_{\mathcal{G}}(c,0)$ when $c_2=0$.

\begin{theorem}\label{th:c2=0}
	Let $c=(c_1,0,c_3)\in\BF^3$.
	\begin{enumerate}[(i)]
		\item If $c_3=0$, then $L_{\mathcal{G}}(c,0)$ is isomorphic
		to $L_{\Vir}(c_1,0)$
		as $\BN$-graded vertex algebras.
		
		\item If $c_3\ne0$, then $L_{\mathcal{G}}(c,0)$ is isomorphic to
		$L_{\widehat{\fh}}(c_3,0)\otimes L_{\Vir}(c_1-1,0)$
		as $\BN$-graded vertex algebras.
	\end{enumerate}
\end{theorem}

\begin{proof}
(i) By \eqref{eq:inv}, \eqref{eq:Hm-Hn-brac},
and $H(1)\1=0$, we have
\begin{equation}\label{eq:v-H-Bili=0}
	(v_H,v_H)=-(\1,H(1)H(-1)\1)=-c_3=0.
\end{equation}
Since $V_{\mathcal{G}}(c,0)_{(1)}=\BF v_H$,
\eqref{eq:orth-m-n} and \eqref{eq:v-H-Bili=0}
imply $v_H\in\Rad(\cdot,\cdot)$.
The relations \eqref{eq:H-I;H-J}
and $I(-1)\1=J(-1)\1=0$ give
$I(-1)v_H=-v_I$ and $J(-1)v_H=v_J$.
Since $\Rad(\cdot,\cdot)$ is an ideal,
$v_I,v_J\in\Rad(\cdot,\cdot)$.
Thus $L_{\mathcal{G}}(c,0)$ is generated by
$-v_L+\Rad(\cdot,\cdot)$.

The operators $-L(m)$ on $L_{\mathcal{G}}(c,0)$
satisfy the Virasoro relations of central charge $c_1$,
and $L(m)\1=0$ for $m\ge-1$.
Hence the universal property of $V_{\Vir}(c_1,0)$
(see \cite[Remark~6.1.3]{LL})
gives a $\Vir$-module homomorphism
\begin{equation*}
	\varphi:V_{\Vir}(c_1,0)\to L_{\mathcal{G}}(c,0)
\end{equation*}
such that $\varphi(\1)=\1+\Rad(\cdot,\cdot)$.
Then
\begin{equation}\label{eq:vir-mod-homo}
	\varphi\big(L(-m_1)\cdots L(-m_k)\1\big)
	=
L(-m_1)\cdots L(-m_k)\1
	+\Rad(\cdot,\cdot)
\end{equation}
for $k\in\BN$ and $m_1,\ldots,m_k>1$,
where $L(m)$ denotes the action of $L_m$ on both modules.
Since $\omega=-L(-2)\1$ generates $V_{\Vir}(c_1,0)$,
it follows from \cite[Proposition~5.7.9]{LL} that
$\varphi$ is a vertex algebra homomorphism.
As $-v_{L}+\Rad(\cdot,\cdot)=\varphi(\omega)$ generates $L_{\mathcal{G}}(c,0)$,
$\varphi$ is surjective.

By \eqref{eq:vir-mod-homo},
$\varphi$ preserves the $\BN$-grading.
Hence $\ker\varphi$ is a graded ideal
of $V_{\Vir}(c_1,0)$.
As $L_{\mathcal{G}}(c,0)$
has no nonzero proper graded ideals
by Lemma~\ref{lem:rad-max-ideal},
$\ker\varphi$ is a maximal graded ideal
of $V_{\Vir}(c_1,0)$.
The uniqueness of the maximal graded ideal
of $V_{\Vir}(c_1,0)$
(see \cite[Section~5]{LM2}) now gives
\begin{equation*}
	L_{\mathcal{G}}(c,0)
	\cong V_{\Vir}(c_1,0)/\ker\varphi
	\cong L_{\Vir}(c_1,0)
\end{equation*}
as $\BN$-graded vertex algebras.

(ii)
Let $I_n$ and $J_n$ act trivially on
$V_{\mathcal{L}}(c,0)$ for $n\in\BZ$,
extending its $\mathcal{L}$-module structure to a $\mathcal{G}$-module structure of central charge $c$.
Since $X(m)\1=0$ for $m\ge-1$ and $X=L,I,J$, 
and $H(m)\1=0$ for $m\ge0$,
the universal property of $V_{\mathcal{G}}(c,0)$ gives a $\mathcal{G}$-module homomorphism
\begin{equation*}
\psi:V_{\mathcal{G}}(c,0)\to V_{\mathcal{L}}(c,0)
\end{equation*}
such that $\psi(\1)=\1$.
As in (i), \cite[Proposition~5.7.9]{LL} shows that $\psi$ is a vertex algebra homomorphism.
Since $\psi(v_L)=L(-2)\1$ and $\psi(v_H)=H(-1)\1$ generate $V_{\mathcal L}(c,0)$, $\psi$ is surjective.
The PBW theorem shows that $\psi$ preserves the $\BN$-grading.

By Lemma \ref{lem:HV-tensor} and the canonical quotient maps,
$\psi$ gives a surjective homomorphism of $\BN$-graded vertex algebras
\begin{equation*}
	\bar{\psi}:V_{\mathcal{G}}(c,0)\to
	L_{\widehat{\fh}}(c_3,0)\otimes
	L_{\Vir}(c_1-1,0).
\end{equation*}
By \cite[Lemma~2.1(2)]{DLR2},
$L_{\widehat{\fh}}(c_3,0)\otimes L_{\Vir}(c_1-1,0)$ has no nonzero proper graded ideals.
Hence $\ker\bar{\psi}$ is a maximal proper graded ideal of $V_{\mathcal{G}}(c,0)$.
Lemma \ref{lem:rad-max-ideal} therefore gives
$\ker\bar{\psi}=\Rad(\cdot,\cdot)$.
Consequently,
\begin{equation*}
		L_{\mathcal{G}}(c,0)
		\cong V_{\mathcal{G}}(c,0)/\ker\bar{\psi}
		\cong L_{\widehat{\fh}}(c_3,0)\otimes
		L_{\Vir}(c_1-1,0)
\end{equation*}
as $\BN$-graded vertex algebras.
\end{proof}

Let $L_{\Vir}(\tfrac12,h)$ denote the irreducible
	graded highest weight $\Vir$-module of central charge $\tfrac12$
	and highest weight $h$ (see \cite[Lemma~3.1]{DR2015}).	
Set
\begin{equation*}
		S_p:=
		\begin{cases}
			\{0,\tfrac12,\tfrac1{16}\}, & \text{if }p\ne7,\\
			\{0,\tfrac12\}, & \text{if }p=7.
		\end{cases}
\end{equation*}
By \cite[Theorems~4.6, 5.2]{DR2015},
	the irreducible $\BN$-graded $L_{\Vir}(\tfrac12,0)$-modules are,
	up to equivalence, precisely $L_{\Vir}(\tfrac12,h)$ with
	$h\in S_{p}$.
	
Combining Theorem~\ref{th:c2=0},
\cite[Theorem~5.8]{LM},
\cite[Theorems~4.6, 4.9, 5.2]{DR2015},
and \cite[Lemma~2.1(2), (3)]{DLR2},
we obtain the following result.

\begin{corollary}\label{coro:c2=0}
For each choice of $c\in \BF^3$ below,
$L_{\mathcal{G}}(c,0)$ has,
up to equivalence, precisely the following
	irreducible $\BN$-graded modules:
\begin{enumerate}[(i)]
	\item For $c=(\frac12,0,0)$:
	$L_{\Vir}(\tfrac12,h)$, $h\in S_p$.
	
	\item 
		For $c=(\frac32,0,c_3)$ with $c_3\ne0$:
	$L_{\widehat{\fh}}(c_3,0)\otimes
	L_{\Vir}(\tfrac12,h)$, $h\in S_p$.
\end{enumerate}
Moreover, 
for $c$ as in (i) or (ii),
every $\BN$-graded $L_{\mathcal{G}}(c,0)$-module
is completely reducible if $p\ne7$.
\end{corollary}

\end{document}